\documentclass[12pt]{article}
\usepackage{amsmath,amsthm,amsfonts,amscd,amssymb,eucal}
\begin{document}
%%%%SIMBOLI%%%%%%%%%%%%%%%%%%%%%%%%%%%%%%%%%%%%%%%%%%
\def\RR{{\mathbb R}}
\def\CC{{\mathbb C}}
\def\NN{{\mathbb N}}
\def\ZZ{{\mathbb Z}}
\def\s1{{S^1}}
\def\diff{{{\rm Diff}^+ (S^1)}}
\def\vect{{{\rm Vect}(S^1)}}
\def\sl2{{{\rm SL}(2,\RR)}}
\def\psl2{{{\rm PSL}(2,\RR)}}
\def\mob{{\rm Mob}}
\def\u1{{{\rm U}(1)}}
\def\su2{{{\rm SU}(2)}}
\def\suN{{{\rm SU}(N)}}
\def\so3{{{\rm SO}(3)}}
\def\gk {{G_k}}
\def\A{{\cal A}}
\def\B{{\cal B}}
\def\C{{\cal C}}
\def\D{{\cal D}}
\def\F{{\cal F}}
\def\H{{\cal H}}
\def\I{{\cal I}}
\def\K{{\cal K}}
\def\k{{\rm K}}
\def\M{{\cal M}}
\def\N{{\cal N}}
\def\O{{\cal O}}
\def\P{{\cal P}}
\def\R{{\cal R}}
\def\S{{\cal S}}
\def\T{{\cal T}}
\def\U{{\cal U}}
\def\V{{\cal V}}
\def\W{{\cal W}}
\def\G{{\bf G}}

%%%%%%%%%%%%%Teoremi%%%%%%%%%%%%%%%%%%%
\newtheorem{theorem}{Theorem}[section]
\newtheorem{definition}[theorem]{Definition}
\newtheorem{corollary}[theorem]{Corollary}
\newtheorem{proposition}[theorem]{Proposition}
\newtheorem{lemma}[theorem]{Lemma}
\newtheorem{remark}[theorem]{Remark}
%%%%%%%%%%%%%%%%%%%%%%%%%%%%%%%%%%%%%%%%%%%%%%%%%%%%
\title{{\huge {\bf On the Representation Theory of Virasoro Nets}}}

\author{SEBASTIANO CARPI\footnote{Supported in part by the Italian 
MIUR and GNAMPA-INDAM.}
\\
Dipartimento di Scienze \\
Universit\`a ``G. d'Annunzio" di Chieti-Pescara\\
Viale Pindaro 42, I-65127 Pescara, Italy\\
e-mail carpi@sci.unich.it}
\date{}
\maketitle
\renewcommand{\sectionmark}[1]{}
\begin{abstract} 
We discuss various aspects of the representation theory of the local nets 
of von Neumann algebras on the circle associated with positive 
energy representations
of the Virasoro algebra (Virasoro nets). 
In particular we classify the local extensions of the $c=1$ Virasoro net 
for which the restriction of the vacuum representation 
to the Virasoro subnet is a direct sum of irreducible subrepresentations with 
finite statistical dimension (local extensions of compact type).  
Moreover  we prove that if the central charge $c$ is in 
a certain subset of $(1,\infty)$, including $[2,\infty )$, and 
$h \geq (c-1)/24$, the irreducible representation with lowest weight 
$h$ of the corresponding Virasoro net has infinite statistical 
dimension. As a consequence we show that if the central charge $c$ is in 
the above set and satisfies $c\leq 25$ then the corresponding Virasoro 
net has no proper local extensions of compact type.
\end{abstract}

\section{Introduction}
The idea that  the formulation of relativistic quantum physics in terms of 
local nets of von Neumann algebras (see e.g. \cite{Haag}) provides a natural 
framework for the classification of 
two-dimensional conformal field theories was already present in the 
late eighties in the work of Buchholz, Mack and Todorov \cite{BMT}. 
As an illustration of this idea these authors classified the 
local conformal nets over $\s1$ (compactified light ray) whose common 
``germ" is the $\u1$ chiral current algebra, namely the local nets 
extending 
the one generated by a $\u1$ current. In the same paper they suggested 
the more general (and ambitious) classification  program of  
conformal field theories $\s1$ whose ``germ" is the Virasoro algebra 
Vir. 
In other words, they proposed to classify the local extensions of the Virasoro 
nets, i.e. the local nets of von Neumann algebras 
on $\s1$ which are generated by the positive energy unitary 
irreducible representations 
with lowest weight $0$ (vacuum representations) of Vir or equivalently 
(see \cite{mack}) by a chiral energy-momentum tensor $T(z)$, cf. \cite{BS-M}.
Since the equivalence class of a Virasoro net is completely determined by 
the value of the central charge $c$ in the corresponding representation of 
Vir, one has to classify the local extensions of a family of nets 
labelled by a positive real number $c$ and this is a clearly well defined 
problem which in fact turns out to be equivalent to the one of 
classifying diffeomorphism covariant nets on the circle.

In a recent remarkable paper \cite{KL} Kawahigashi and Longo have been 
able to solve the above problem for all the Virasoro nets with $c<1$ and 
subsequently they used this result to classify all local conformal nets 
on the two-dimensional Minkowski space-time, having parity symmetry and 
central charge less then 1 \cite{KL2}. The extension of their results 
in the $c\geq 1$ region appears to be a very important and 
difficult challenge. 

In the transition from $c<1$ to $c\geq 1$ two drastic differences 
are immediately evident.
 The first is that the Virasoro nets with $c\geq 1$ 
are all known to be non-rational. The second is that they are all expected 
to have irreducible sectors with infinite statistical dimension 
\cite{Reh}, a fact that has been proved by the present author in the case 
$c=1$ \cite{Car02}. Rationality  and the absence of irreducible sectors 
with infinite statistical dimension play a fundamental role in the 
classification of $c<1$ conformal nets and hence, some of the main ideas 
in \cite{KL} do not apply for the remaining values of the central charge. 

The purpose of this paper is give some new insight in the understanding 
of the representation theory of $c\geq 1$ Virasoro nets and their local 
extensions 
with the above problems in mind, especially concerning the role of 
infinite statistical dimension. 
Our main results are the classification of local extensions of compact 
type (see Definition \ref{defcompacttype}) for the 
Virasoro net with $c=1$ (Theorem \ref{classc=1})
\footnote{An equivalent result as been independently 
obtained by F. Xu \cite{xu 2003}} 
and the proof that if the central 
charge $c$ is in a certain subset of $(1,\infty)$ containing $[2,\infty)$, 
then  the irreducible positive energy representations with lowest weight 
$h\geq (c-1)/24$ of the corresponding Virasoro net have infinite 
statistical dimension (Theorem \ref{infdim} ), a fact that it is expected 
to hold for every $c>1$ and $h>0$ \cite{Reh}. As a consequence of the 
latter result we show that if $c\in [2,25]$ then the corresponding 
Virasoro net as no proper local extensions of compact type 
(Theorem \ref{cmaxnon-comp}). As in the $c=1$ case \cite{Car02}, we use 
oscillator (Fock) representations of Vir to obtain the result on infinite 
statistical dimension but the argument we have found for $c>1$ is more 
intricate and relies in part on  recent results of S. K\"oster 
\cite{koester03b}.    

Besides these main results we also provide the proof of some relevant 
properties of the Virasoro nets  which seem to appear at  most 
implicitly 
in the literature, like the fact that every irreducible positive energy 
representation of a Virasoro net on a separable Hilbert space comes from a 
representation of the Virasoro algebra (Prop. \ref{allvirrep}) or the 
fact that every local extension of a Virasoro net is diffeomorphism 
covariant (Prop. \ref{extvir}).

\section{Preliminaries}
\label{preliminaries}
\subsection{Conformal nets on $\s1$ and diffeomorphism covariance}
\label{prel1}
Let $\I$ be the set of nonempty, nondense,  open  intervals of the
unit circle $\s1 =\{z\in \CC :|z|=1 \}$.

A {\bf conformal net on $\s1$} is a family
$\A=\{ \A(I) : I \in \I \}$
of von Neumann algebras, acting on a infinite-dimensional separable
complex Hilbert space ${\H}_\A$, satisfying the following properties: 
\begin{itemize}
\item[(i)] {\it Isotony.}
\begin{equation}
\A(I_1) \subset \A (I_2),\;{\rm if}\; I_1 \subset I_2,
\;
I_1, I_2 \in \I.
\end{equation}

\item[(ii)] {\it Locality.}
\begin{equation}
\A(I_1) \subset \A(I_2)',\;{\rm if}\;
I_1 \cap I_2 =\emptyset ,\;I_1, I_2 \in \I.
\end{equation}

\item[(iii)] {\it M\"obius covariance.} There exists  a strongly
continuous unitary representation $U$ of $\psl2$ in ${ \cal H}$ such that
\begin{equation}
U(g){\cal A}(I)U(g)^{-1}={\cal A}(gI),\; I\in \I,\; g\in \psl2,
\end{equation}
where $\psl2$ acts on $\s1$ by M\"obius transformations (cf. Appendix 
\ref{appendix}). 
 
\item[(iv)] {\it Positivity of the energy.} The representation $U$ has 
positive energy, namely the conformal Hamiltonian $L_{0}$,
which generates the restriction of $U$
to the one-parameter subgroup of rotations $r(\vartheta)$,
has nonnegative spectrum.

\item[(v)] {\it Existence and uniqueness of the vacuum.} There exists a 
unique (up to aphase)  $U$-invariant unit vector $\Omega\in \H_\A$.

\item[(vi)] {\it Cyclicity of the vacuum.} ${\Omega}$ is cyclic for the 
algebra $\A(\s1):=\bigvee_{I\in \I} \A (I)$ 
\end{itemize} 

Some consequences of the axioms are \cite{FrG,GuLo96,FJ}:
\begin{itemize}
\item[(vii)] {\it Reeh-Schlieder property.} For every $I\in \I$,
$\Omega$ is cyclic and separating for $\A(I)$.

\item[(viii)] {\it Bisognano-Wichmann property.} If ${\Delta}_{I}$ is the
modular operator associated to ${\A}(I)$ and $\Omega$  then
\begin{equation}
{\Delta}_{I}^{it}=U(\Lambda_I ({(2\pi)t}),
\end{equation}  
where $\Lambda_I$ is the one parameter subgroup of $\psl2$ of special
conformal transformations preserving $I$.

\item[(ix)] {\it Haag duality.} For every $I \in \I$
\begin{equation}
\A(I)'= \A(I^{c}),
\end{equation}
where $I^{c}$ denotes the interior of $S^{1}\backslash I$.

\item[(x)] {\it Factoriality.} The algebras ${\A}(I)$ are type
${\rm III}_{1}$ factors.

\item[(xi)] {\it Irreducibility.} $\A(\s1) = {\rm B}(\H_\A),$ where 
${\rm B}(\H_\A)$ denotes the algebra of all bounded linear operators on 
$\H_\A$.
\item[(xii)]{\it Additivity.} If $\S\subset \I$ is a covering of the 
interval $I$ then  
\begin{equation}
\A(I)\subset \bigvee_{J\in \S}\A(J).
\end{equation}
\end{itemize}
Furthermore, it follows easily from the strong continuity of the 
representation $U$ that conformal nets are {\bf continuous from 
outside}, namely 
\begin{equation}
\A(I) = \bigcap_{J\supset \overline{I}}\A(J).
\end{equation}

A conformal net $\A$ is said to be  {\bf split} if given two 
intervals
$I_1,I_2\in\I$ such that the closure $\overline{I_1}$ of $I_1$ is 
contained in $I_2$, there exists
a type I factor $\N(I_1,I_2)$ such that
\begin{equation}
\A(I_1)\subset \N(I_1,I_2) \subset \A(I_2).
\end{equation}
If ${\rm Tr}(t^{L_0})<\infty$ for every $t\in (0,1)$ then $\A$ is
split \cite[Theorem 3.2]{D'ALR}.

$\A$ is said to be {\bf strongly additive} if for every $I,I_1,I_2\in\I$ 
with $I_1,I_2$ obtained by removing a point from $I$ we have
\begin{equation}
\A(I_1)\vee\A(I_2)=\A(I).
\end{equation}
It is often convenient to identify $\s1/\{-1\}$ with the real line $\RR$. 
With this identification, the family of nonempty open bounded 
intervals of $\RR$ corresponds to the family 
$\I_0=\{I\in \I: -1 \notin \overline{I}\}.$
The restriction $\A_0$ of a conformal net $\A$ to $\I_0$ can be considered 
as a net on $\RR$. Moreover, since $\I_0$ is directed under inclusion,
one can define the quasi-local $C^*$-algebra 
(still denoted $\A_0$) $(\bigcup_{I\in \I_0}\A(I))^{-||\cdot ||})$
as $C^*$-inductive limit of the local von Neumann algebras $\A(I)$, 
$I\in \I_0$.
\medskip 

We now briefly discuss diffeomorphism covariance. Let $\diff$ the group 
of orientation preserving diffeomorphisms of the circle. It is an infinite 
dimensional Lie group modelled on the real topological vector space 
$\vect$ of smooth real vectors fields on $\s1$ with the usual $C^\infty$ 
topology \cite[Sect.6]{Milnor}. Its Lie algebra coincides with $\vect$ 
with the bracket given by the negative of the usual brackets of vector 
fields. Hence if $g(z), f(z)$, $z=e^{i\vartheta}$, are real 
valued functions in $C^\infty(\s1)$ then 
\begin{equation}
[g(e^{i\vartheta})\frac{d}{d\vartheta}, 
f(e^{i\vartheta})\frac{d}{d\vartheta}]= 
(\frac{d}{d\vartheta}g(e^{i\vartheta}))f(e^{i\vartheta})-
(\frac{d}{d\vartheta}f(e^{i\vartheta}))g(e^{i\vartheta}).
\end{equation} 
In this paper we shall often identify the vector field 
$g(e^{i\vartheta})\frac{d}{d\vartheta}\in \vect$ with the corresponding 
real function $g(z) \in C^\infty(\s1)$. 

Following \cite{KL} for every $I\in \I$ we shall denote by 
${\rm Diff}(I)$ the subgroup of $\diff$ whose elements are the 
diffeomorphisms of the circle which act as the identity on $I$. Note that 
${\rm Diff}(I)$ does not coincide with the group of diffeomorphisms of 
the open interval $I$, as the notation might erroneously suggest.

By a strongly continuous projective unitary representation 
$V$ of $\diff$ on a Hilbert space we shall always mean a strongly 
continuous homomorphism of $\diff$ into the quotient 
${\rm U(\H)}/\mathbb T$ 
of the unitary group of $\H$ by the circle subgroup 
$\mathbb T$. The restriction of the representation $V$ to the M\"obius 
subgroup of $\diff$ always lifts to a unique strongly continuous unitary 
representation 
$U$ of the universal covering group $\widetilde{\psl2}$ of $\psl2$. 
We shall say 
that $V$ extends $U$ and that $V$ is a positive energy representation if 
$U$ is a positive energy representation of $\widetilde{\psl2}$, namely 
if the corresponding conformal Hamiltonian $L_0$ , which generates the 
restriction of $U$ to the lifting $\tilde{r}(\vartheta)$ of the
one-parameter subgroup $r(\vartheta)$ of rotations, has nonnegative
spectrum. Note that although for $\gamma \in \diff$, $V(\gamma)$ is 
defined only up to a phase as an operator on $\H$, expressions 
like $V(\gamma)TV(\gamma)^*$ for $T\in {\rm B}(\H)$ or $V(\gamma) \in \M$ 
for a (complex) subspace $\M \subset {\rm B}(\H)$ are unambiguous and will 
be used in the following.  

We shall say that a conformal net on $\s1$ is {\bf diffeomorphism
covariant} if there is a strongly continuous projective unitary 
representation $V$ of $\diff$ on $\H_\A$ extending $U\circ q$ ( where 
$U$ is the original unitary representation  of $\psl2$ making $\A$ 
M\"obius covariant and 
 $q: \widetilde{\psl2} \mapsto \psl2$ denotes the covering map)
 and such that, for every $I\in \I$ 
\begin{eqnarray}
V(\gamma)\A(I)V(\gamma)^*= \A(\gamma I),\;\gamma \in \diff, \\
V(\gamma)AV(\gamma)^*=A, \;  \gamma \in {\rm Diff}(I), A\in \A (I^c).
\end{eqnarray}

\subsection{Representations of conformal nets}
\label{prel2}

A {\bf representation} of a conformal net $\A$ is a family 
$\pi =\{ \pi_I : \; I\in \I \}$ where $\pi_I$ is a (unital) representation 
of $\A(I)$ on a fixed Hilbert space $\H_\pi$, such that 

\begin{equation}
\pi_J |_{\A(I)} =\pi_I \;\; {\rm if} \;I\subset J,\; I, J\in\I.
\end{equation}

Irreducibility, direct sums and unitary equivalence of representations
of a conformal net can be defined in a natural way, see \cite{FrG,GuLo96}. 

If $\H_\pi$ is separable then, since the local von Neumann algebras 
$\A(I), I\in \I$, are factors, $\pi$ is automatically 
locally normal, namely $\pi_I$ is normal for each $I\in \I$, see 
\cite{takesaki}. Hence, $\pi_I(\A(I))$ is a type 
${\rm III}_1$ factor.
The unitary equivalence class a representation $\pi$ on a separable
Hilbert space is
called a {\bf sector} and denoted $[\pi]$. If $\pi$ is irreducible then we
say that $[\pi]$ is an irreducible sector (also called superselection
sector). The defining representation $\pi_0$ of a conformal net $\A$ on
the Hilbert space $\H_\A$ is called the {\bf vacuum representation}. The 
corresponding sector is called the vacuum sector and $\H_\A$ is said to 
be the vacuum Hilbert space of $\A$.

A representation $\pi$ is said to be {\bf covariant} if there is a 
strongly continuous unitary representation $U_\pi$ of $\widetilde{\psl2}$ 
on $\H_\pi$ such that 

\begin{equation}
\label{covariance}
{\rm Ad} U_\pi (g)\circ\pi_I=\pi_{q(g) I}\circ{\rm Ad} U(q(g)),\; g\in
\widetilde{\psl2}, I\in \I.
\end{equation}

If $U_\pi$ can be chosen 
to be a positive energy representation, 
then $\pi$ is said to be {\bf covariant with positive energy}. 
In this case one can always choose $U_\pi$ to be with
positive energy and inner, namely such that 
\begin{equation}
\label{innercov}
U_\pi(\widetilde{\psl2})\subset \pi(\A)'':= \bigvee_{I\in \I}\pi_I(\A(I)), 
\end{equation} 
and this choice is unique, see \cite{koester02} and (the proof of) 
\cite[Lemma 5.14]{BCL}. 

Given a covariant representation $\pi$ of $\A$ on a separable Hilbert
space 
$\H_\pi$ one has the (isomorphic) type III subfactors 
$\pi_I(\A(I))\subset \pi_{I^c}(\A(I^c))',\; I\in \I$ \cite{FrG}. Hence the 
corresponding (minimal) index $[ \pi_{I^c}(\A(I^c))' : \pi_I(\A(I)) ]$ 
\cite{jones,kosaki,longo 89-90}
is independent of $I\in \I$ and the {\bf statistical dimension} $d(\pi)$ 
of $\pi$ is given by 
\begin{equation}
\label{def.dim.}
d(\pi )= [ \pi_{I^c}(\A(I^c))' : \pi_I(\A(I))]^{\frac{1}{2}}.
\end{equation}  
\medskip

A representation $\rho$ of a conformal net $\A$ on its vacuum Hilbert 
space $\H_\A$ is said to be {\bf localized} in an interval $I_0\in \I$
if $\rho_{I_0^c}$ is the identity representation. As a consequence of Haag 
duality if a representation $\rho$ is localized in $I_0$ and $I\in \I$ 
contains $I_0$ then $\rho_I$ is an endomorphism of $\A(I)$ whose index 
is the square of the statistical dimension of the representation $\rho$.
Moreover, for every interval $I_0\in \I$ and every representation 
$\pi$ of $\A$ on a separable Hilbert space one can find a representation 
$\rho$, localized in $I_0$ and unitarily equivalent to $\pi$, see 
\cite{FrG,GuLo96}. 
The restriction to $\A_0$ of a representation $\rho$ localized in some 
$I_0\in \I_0$ is called a {\bf DHR endomorphism} and in fact yields
an endomorphism of the quasi-local $C^*$-algebra $\A_0$, see 
\cite[Sect. 3]{longo 2003}. The set of DHR endomorphisms is a semigroup 
(under composition) and it has a natural (DHR) unitary braiding, see 
\cite{FRS2,FrG}. As usual we shall denote $\epsilon(\rho,\sigma)$ the 
unitary braiding operator associated to the DHR endomorphisms $\rho$ and 
$\sigma$.

\subsection{Subsystems}
\label{prel3}
A {\bf conformal subsystem}  (or subnet) of a conformal net
${\cal A}$ is a family 
${\B} =\{ {\B}(I) :  I\in \I \}$ of nontrivial von Neumann
algebras acting on ${\cal H}_\A$ such that:
\begin{eqnarray}
{\B}(I)\subset {\A}(I) \; I\in \I;\\
U(g){\B}(I)U(g)^{-1}={\B}(g I)
\; I \in \I,\; g\in \psl2;\\  
{\B}(I_1)\subset {\B}(I_2)
\;{\rm if}\;I_1 \subset I_1, \;\;I_1, I_2 \in \I.
\end{eqnarray}

We shall use the notation $\B\subset \A$ for conformal subsystems. 
Note that $\B$ is not in general a conformal net since $\Omega$ is not 
cyclic for the algebra $\B(\s1):=\bigvee_{I\in \I} \B (I)$, unless 
$\B=\A$.
However one gets
a conformal net $\widehat{\B}$ restricting the algebras $\B(I),\;I\in\I,$
and of the representation $U$ to the closure $\H_\B$ of $\B(\s1)\Omega$.
Since the map $b\in \B(I) \mapsto b |_{\H_\B} \in \widehat{\B}(I) $
is an isomorphism for every $I\in \I$, because of the Reeh-Schlieder
property, as usual, we shall often use the symbol $\B$ 
instead of $\widehat{\B}$, specifying, if necessary, when $\B$ acts on
$\H_\A$ or on $\H_\B$.  

Let $\pi$ be the defining representation of the conformal net 
$\B\subset \A$ on the Hilbert space $\H_\A$ (i.e. the restriction to $\B$ 
of the vacuum representation of $\A$). Because of the separability of 
$\H_\A$, for every $I_0\in \I$ we can find a representation $\theta$ on the 
vacuum Hilbert space $\H_\B$ of $\B$, which is  unitarily equivalent to $\pi$
and is localized in $I_0$. Then if $I_0 \subset I\in \I$, $\theta_I$ is a 
dual canonical endomorphism for the subfactor $\B(I)\subset \A(I)$, namely 
there is a canonical endomorphism 
(in the sense of \cite{longo 89-90}) for the latter whose restriction to 
$\B(I)$ coincides with $\theta_I$, see \cite[Proposition 3.4]{LR} and
\cite[Sect. 3.3]{longo 2003}.

\subsection{Virasoro nets and their representations}
\label{virsubsec}

Let Vir denote the Virasoro algebra i.e. the complex Lie algebra 
spanned by $L_n, n\in \ZZ$ and a central element $\kappa$ with relations
\begin{equation}
\label{lievir}
[L_n,L_m]=(n-m)L_{n+m} +\delta_{n+m,0} \frac{n^3-n}{12} \kappa.
\end{equation}
We shall denote $L(c,h)$ the unique positive energy irreducible unitary 
representation of Vir with lowest weight $h$ and central charge $c$ 
(see e.g. \cite{DMS,KaRa}). 
The conformal Hamiltonian $L_0$ is diagonalizable on the corresponding 
Vir-module (still denoted $L(c,h)$) with spectrum $h+\NN_0$ and  
the central element $\kappa$ acts as multiplication by the real number 
$c$. Positivity of the energy implies $h\geq 0$ and
unitarity (or hermiticity) means that there is a positive definite
sesquilinear form $(\cdot,\cdot)$ on $L(c,h)$ such that 
\begin{equation}
(\xi,L_n \psi)=(L_{-n}\xi,\psi),
\end{equation}
for $\xi, \psi \in L(c,h)$, $n\in\ZZ$. 

The above conditions  
give restrictions on the values of the pair $(c,h)$. 
In fact either $c\geq 1$ and $h\geq 0$ or 
we have a pair $(c(m),h_{p,q}(m))$, $m\in \NN$, where 
\begin{equation}
\label{discretec}
c(m) = 1- \frac{6}{(m+2)(m+3)}
\end{equation}
and 
\begin{equation}
h_{p,q}(m)=\frac{((m+3)p-(m+2)q)^2-1}{4(m+2)(m+3)},
\end{equation}
$p=1,...,m+1$, $q=1,...,p$, (discrete series representations).
For later convenience we shall denote $D\subset [\frac{1}{2},1)$ the set 
of discrete values of the central charge in Eq. (\ref{discretec}). 
Accordingly the set of allowed values of the central charge is 
$D\cup [1,\infty)$.  

Now let $\H(c,h)$ be the Hilbert space completion of the module 
$L(c,h)$. Then the Virasoro algebra acts on $\H(c,h)$ by unbounded 
operators on the common invariant domain $L(c,h) \subset \H(c,h)$ which 
can in fact be identified with the subspace $\H^{fin}(c,h)$ of finite 
energy vectors i.e. the linear span of the eigenvectors of the conformal 
Hamiltonian. 
The (chiral) energy-momentum tensor $T_{(c,h)}(z)$, $z=e^{i\vartheta}\in 
\s1$ associated to $L(c,h)$, is defined by the formal power 
series
\begin{equation}
T_{(c,h)}(z)=\sum_{n\in \ZZ} L_n z^{-n-2}
\end{equation}

For a function on $\s1$, $\vartheta \mapsto f(e^{i\vartheta})$ with finite 
Fourier series (trigonometric polynomial), the 
operator
\begin{equation}
T_{(c,h)}(f)= \sum_{n\in \ZZ}L_n f_{n},
\end{equation}
where 
\begin{equation}
f_n=\int_0^{2\pi}\frac{d\vartheta}{2\pi}e^{-in\vartheta}f(e^{i\vartheta}),
\end{equation}
belongs to Vir and hence is well defined on $\H^{fin}(c,h)$ and leave it 
invariant. 
Also the following (formal) notation is used
\begin{equation}
T_{(c,h)}(f)= \oint_\s1 \frac{zdz}{2\pi i}T_{(c,h)}(z)f(z)
=\int_0^{2\pi}\frac{d\vartheta}{2\pi} 
T_{(c,h)}(e^{i\vartheta})e^{i2\vartheta}f(e^{i\vartheta}).
\end{equation}

The Virasoro net $\A_{({\rm Vir},c)}$ can be defined as in \cite{BS-M} 
as the net generated by the energy-momentum tensor $T_c(z) :=T_{(c,0)}(z)$ 
in the 
representation of lowest weight $0$ on 
$\H(c,0)=: \H_{\A_{({\rm Vir},c)}}$. 
First of all one can show that 
the map $f \mapsto T_c(f)$ extends (uniquely) to an operator valued 
distribution
(Wightman field) on the invariant domain $C^\infty(L_0)$, the subspace of 
smooth vectors for $L_0$. Moreover the linear energy-bounds established in
 \cite{BS-M} (also cf. \cite{GoWa}) imply that for every smooth 
real valued function $f$, $T_c(f)$ is essentially 
self-adjoint (on any core for $L_0$) and that $e^{iT_c(f_1)}$
commutes with $e^{iT_c(f_2)}$ when the real smooth functions $f_1$ and 
$f_2$ have disjoint supports 
(in fact these properties also hold in the representations with $h>0$). 
It follows that the net of von Neumann algebras defined by    
\begin{equation}
\A_{({\rm Vir},c)} (I)= \{{\rm e}^{iT_c(f)}: f\in C^\infty(\s1),\; 
{\rm real},\; {\rm supp}\; f\subset I\}'',\;\;I\in\I.
\end{equation}
is local and in fact one can verify all the other axioms of a conformal 
net. In particular the representation $U$ of $\psl2$ is obtained by 
integrating the self-adjoint part of the (complex) Lie subalgebra of 
Vir spanned by $L_{-1}, L_0, L_1$ and the vacuum vector $\Omega$ is 
the (normalized) lowest weight vector in $L(c,0)$. 

An alternative construction is obtained by integrating the representations
$L(c,0)$ of Vir to the corresponding projective unitary representations 
of $\diff $. In fact as shown by Goodman and Wallach \cite{GoWa}
(cf. also \cite{Tol99}), for each allowed pair $(c,h)$
 there is a unique strongly continuous projective unitary representation 
$V_{(c,h)}$ of $\diff$ on $\H(c,h)$ satisfying 
\begin{equation}   
\label{diffvir}
V_{(c,h)}(\exp(f)) = p( e^{iT_{(c,h)}(f)} )
\end{equation}
for every real smooth function (vector field) $f$ on $\s1$. Here 
$\exp(f)\in \diff$ denotes the exponential of the vector field $f$, namely
$t \mapsto \exp(tf)$ is the unique one-parameter group of diffeomorphisms 
generated by $f$, and 
$p: {\rm U}(\H(c,h)) \mapsto {\rm U}(\H(c,h))/{\mathbb T}$ denotes 
the quotient map.
Then the net $\A_{({\rm Vir},c)}$ can be defined by 
\begin{equation}\A_{({\rm Vir},c)}(I) = \{V_{(c,0)}(\gamma) : 
\gamma \in {\rm Diff}(I) \}'', 
\end{equation} 
$I\in \I$. 

The two definitions are equivalent because the group generated by 
the exponentials of smooth vector fields with support in $I\in \I$ is 
dense in ${\rm Diff}(I)$, see \cite[Sect. V.2]{loke}. From the second one 
the diffeomorphism covariance of the Virasoro nets is explicit.

As a consequence of the finiteness of the (vacuum) Virasoro characters 
$\chi(t):= {\rm Tr}(t^{L_0})$ for every $t\in (0,1)$ the Virasoro nets are
split for every allowed value of the central charge. 
For $c\leq 1$ the Virasoro nets are strongly additive \cite{KL,xu 2003} 
while for $c>1$ they are not \cite{BS-M}.
\medskip 

We now discuss some properties of the representation theory of the 
Virasoro nets that we shall need in the following. Let $\H(c,h)$ be the 
Hilbert space completion of Vir module $L(c,h)$ as at the beginning of 
this subsection and let $T_{(c,h)}(z)$ be the corresponding 
energy-momentum tensor. 
A representation of $\A_{({\rm Vir},c)}$ on $\H(c,h)$ will 
be denoted $\pi^c_h$ if for every $I\in \I$ and every real smooth real 
function $f$ on $\s1$ with support in $I$, the following hold 
\begin{equation}
\label{pich} 
{\pi^c_h}_I (e^{iT_c(f)}) = e^{iT_{(c,h)}(f)}.
\end{equation}
 
It is immediate to verify that if  a representation satisfying Eq. 
(\ref{pich}) exists, then it is unique. More complicate is to 
demonstrate the existence of such representations. Of course the vacuum 
representation $\pi^c_0$ exists for every allowed value of $c$ i.e. 
for each $c\in D \cup [1,\infty)$. 
If $c<1$ and $h$ is a corresponding allowed value of the lowest 
weight then the representation $\pi^c_h$ exists as a consequence of the
Goddard, Kent, Olive coset construction \cite{GKO} and the local 
equivalence of positive energy representations of the loop groups 
$L\su2$ at fixed 
level \cite{FrG,Was A}, cf. \cite[V.3.3.2]{loke} and \cite[Sect. 3]{KL}.  
If $c\geq 1$ the existence of $\pi^c_h$ has been proved by D. Buchholz and 
H. Schulz-Mirbach for every $h\geq (c-1)/24$. Finally if 
$c\in (D+1)\cup [2,\infty)$, then $c-1$ is an allowed value of the central 
charge. Then using the embedding 
$$\A_{({\rm Vir},c)}\subset \A_{({\rm Vir},c-1)} \otimes \A_{({\rm Vir},1)}$$
one can easily construct, for every $h\geq 0$, the representation 
$\pi^c_h$ as a subrepresentation of the restriction to 
$\A_{({\rm Vir},c)}$ of $\pi^{c-1}_0 \otimes \pi^1_h$. \footnote{I learned 
this argument in an unpublished manuscript of D. Buchholz \cite{Buch90}. } 
As far as we know 
the existence of the representation $\pi^c_h$ for the remaining 
allowed pairs $(c,h)$ is still an important open problem.

\begin{proposition}
\label{allvirrep}
If $\pi$ is an irreducible covariant representation with positive energy 
of the Virasoro net $\A_{({\rm Vir},c)}$ on a separable Hilbert space 
$\H_\pi$ then it is unitarily equivalent to $\pi^c_h$ for some $h\geq 0$. 
\end{proposition}

\begin{proof} Let $V_{(c,0)}$ be the unique projective unitary 
representation of 
$\diff$ on $\H_{\A_{({\rm Vir},c)}}$ such that Eq. (\ref{diffvir}) holds 
with $h=0$. 
From \cite[Sect. 2]{koester03a} (cf. also \cite[Lemma 3.1]{KL}) we know 
that there is a strongly continuous positive energy projective unitary 
representation $V_\pi$ of $\diff$ on $\H_\pi$ such that 
$p(\pi_I(V_{(c,0)}(\gamma)))=V_\pi(\gamma)$ for each $I\in\I$ and 
$\gamma \in {\rm Diff}(I)$. Then it follows from the irreducibility and 
local normality of  $\pi$ that $V_\pi$ is irreducible. 
As a consequence of Theorem \ref{diffdiff} in the Appendix , there is on 
$\H^{fin}_\pi$ 
a positive energy representation $R_\pi$ of the Virasoro algebra with 
central charge $c' \in D\cup [1,\infty)$, which is unitarily equivalent to 
$L(c',h)$ for some $h\geq0$.
Let 
$$T^\pi(z)=\sum_{n\in \ZZ} L^\pi_n z^{-n-2}$$
be  the corresponding energy-momentum tensor. Then, for every real smooth 
vector field $f$ on $\s1$, we have $$V_\pi(\exp (f)) =p(e^{iT^\pi(f)}).$$ 
It follows that if $I\in \I$ and the support of $f$ is contained in $I$ 
$$\pi_I(e^{iT_c(f)})= e^{i\alpha_I(f)}e^{iT^\pi(f)}, $$ 
where $\alpha_I (f)$ is a real constant. Now, it is fairly easy to check 
that there is a (necessarily unique ) distribution $\alpha$  
such that $\alpha(f)= \alpha_I (f)$ 
for every $I\in \I$ and every real function $f $ with support in contained 
in $I$ and that M\"obius covariance implies that $\alpha = 0$. Hence we 
have the equality 
$$\pi_I(e^{iT_c(f)})= e^{iT^\pi(f)},$$
which implies $c'=c$. The conclusion then follows because 
the representation  $R_\pi$ of Vir is unitarily equivalent to  $L(c,h)$ 
for some $h\geq 0$.
\end{proof}
We conclude this subsection with the following proposition. 
\begin{proposition}
\label{directsum}
Let $\pi$ be a positive energy covariant representation of the Virasoro 
net $\A_{({\rm Vir},c)}$ on a separable Hilbert space and let $U_\pi$ be the 
corresponding unique inner unitary representation of 
$\widetilde{\psl2}$. Assume that $U_\pi(\tilde{r}(2\pi))\in \CC 1$. Then 
the following hold: 
\begin{itemize}
\item[(a)] The representation $\pi$ is a direct sum of irreducible 
covariant positive energy representations.
\item[(b)] There exists a unique strongly continuous projective unitary 
representation $V_\pi$ of $\diff$ on $\H_\pi$ satisfying 
\begin{equation}
\label{vpi0}
p(\pi_I(e^{iT_c(f)}))=V_\pi(\exp(f)),
\end{equation}
for every $I\in \I$ and every real smooth function $f$ with support 
contained in $I$. Moreover, this representation satisfies 
\begin{eqnarray}
\label{vpi}
V_\pi(\gamma) \in \pi_I(\A_{({\rm Vir},c)}(I))\quad \forall I\in \I , 
\forall \gamma \in {\rm Diff}(I) ,\\
\label{vpi2}
V_\pi(q(g)) = p(U_\pi(g)) \quad \forall g\in \widetilde{\psl2}.
\end{eqnarray}
\end{itemize} 
\end{proposition}
\begin{proof} The net $\A_{({\rm Vir},c)}$ has the split property and 
hence, has a consequence of \cite[Proposition 56]{KLM}, $\pi$ has a 
direct integral decomposition 
$$\pi= \int_X^\oplus \pi_\lambda d\mu (\lambda),$$
where, for almost every $\lambda$, $\pi_\lambda$ is an irreducible 
representation of $\A_{({\rm Vir},c)}$ on a separable Hilbert space 
$\H(\lambda)$. 
Since $U_\pi(g)\in \pi(\A_{({\rm Vir},c)})''$ for each $g \in 
\widetilde{\psl2}$ we also have the decomposition 
$$U_\pi(g)=\int_X^\oplus U_\lambda(g)d\mu (\lambda).$$
If $h_\pi$ is the lowest eigenvalue of $L^\pi_0$ we have by assumption 
$U_\pi(\tilde{r}(2\pi))=e^{2\pi ih_\pi}$. Hence $U_\lambda$ is, 
for almost every $\lambda$, a positive energy representation satisfying
$L^\lambda_0 \geq h_\pi$ and $U_\lambda(\tilde{r}(2\pi))=e^{2\pi ih_\pi}$. 
It follows that, for almost every $\lambda$, $\pi_\lambda$ is an 
irreducible 
covariant representation of $\A_{({\rm Vir},c)}$ with positive energy 
which, because of Prop. \ref{allvirrep}, is unitarily equivalent to 
$\pi^c_{h_\pi + n_\lambda}$ for some $n_\lambda \in \NN_0$. Now let 
$X_n=\{\lambda \in X: \pi_\lambda \simeq \pi^c_{h_\pi + n} \}$. Then, it 
follows from \cite[Lemma 60]{KLM}, that $\{X_n: n\in \NN_0\}$ is a family 
of pairwise disjoint measurable subsets of $X$ such that 
$\mu (X\backslash \bigcup_{n\in \NN_0}X_n )=0$. Hence 
$$\pi \simeq \bigoplus_{n\in \NN_0}\int_{X_n}^\oplus \pi_\lambda d\mu 
(\lambda)$$ 
and since $\int_{X_n}^\oplus \pi_\lambda d \mu (\lambda)$ is unitarily 
equivalent to a (possibly zero) multiple of $\pi^c_{h_\pi + n}$, (a) 
follows. 

Now it follows from (a) and Prop.\ref{allvirrep} that on 
the dense subspace 
$C^\infty(L^\pi_0)$ of smooth vectors  for $L^\pi_0$ there is a 
projective representation $\eta$ of the Lie algebra of smooth real vector 
fields on $\s1$ by essentially skew-adjoint operators satisfying 
$e^{\eta (f) }=\pi_I (e^{iT_c(f)})$ if $I\in \I$ and ${\rm supp}g \subset 
I$. Moreover $\eta$ satisfies the assumptions in 
\cite[Theorem 5.2.1]{Tol99} (cf. the proof of 
\cite[Theorem 6.1.1]{Tol99} and the discussion in 
\cite[Appendix]{koester03b}). Hence it can be integrated to a unique 
strongly continuous projective unitary representation of the covering 
group of $\diff$ which, since by assumption 
$e^{2\pi i L^\pi_0}=e^{2\pi i h_\pi}$, factors through $\diff$ giving a 
representation $V_\pi$ satisfying Eqs. (\ref{vpi0}) and (\ref{vpi2}). The 
remaining claim in (b) then follows easily.   
\end{proof}

\section{Local extensions}
\begin{definition}
\label{deflocext}
We define a {\bf local extension} of a conformal net $\A$ to be 
a conformal net $(\B, U,\H_\B)$ together with a conformal 
subsystem $\C \subset \B$ such that the corresponding conformal net 
$\widehat{\C}$ on $\H_\C$ is isomorphic to $\A$ and such that 
\begin{equation}
\label{generating}
U(\psl2) \subset \C(\s1).
\end{equation}
\end{definition}

In agreement with the notation for conformal subsystems, 
since $\A$ and $\widehat{\C}$ are isomorphic, we shall often identify 
$\A$ and $\C$  and accordingly we shall write $\A \subset \B$ instead 
of $\C \subset \B$.

Condition (\ref{generating}) implies that $\C \subset \B$ is a 
{\bf full subsystem}, namely that 
\begin{equation}
\C(\s1)' \cap \B(I) = \CC 1 \quad I\in \I. 
\end{equation} 
It prevents trivial extensions of the type $\A \subset \A \otimes \C$, 
cf. \cite{BMT}. For finite index subsystems condition (\ref{generating}) 
is automatically satisfied and we don't know any example of a full 
conformal subsystem violating it. 
Note that in the literature the term ``local extension" is often used
in a weaker sense (see e.g. \cite{LR}).

A class of examples of local extensions is obtained by considering 
fix points subsystems under compact group actions. 
More precisely given a conformal net $\B$ and a strongly 
compact group $G$ of (vacuum preserving) internal symmetries of $\B$ 
one can define the fixed point  subsystem $\A \equiv \B^G$. This kind 
of construction is paradigmatic in the algebraic approach to quantum 
field theory, see \cite{DHR0,DoRo90}. 

One has $\A(\s1)' = G''$ (cf. \cite[Theorem 3.6]{DoRo90}) and since $U$
and $G$ commute (see \cite[Lemma 2.22]{FrG}), 
condition (\ref{generating}) is satisfied. 
Hence $\B$ is a local extension of $\A$ in the sense of Definition 
\ref{deflocext}. 

If $\pi$ is the identical representation of $\A$ on $\H_\B$ one has 
\begin{equation}
\label{reduction1}
\pi = \bigoplus_{\xi \in \hat{G}} d(\xi)\pi_\xi,
\end{equation}
where $\hat{G}$ is set of equivalence classes of irreducible unitary
representations of $G$,  the $\pi_\xi$ are mutually inequivalent 
irreducible covariant representations of $\A$ 
(with trivial univalence) 
appearing with multiplicity $d(\xi)$ equal to the dimension of the 
representations of G of class $\xi$ and satisfying $d(\pi_\xi)=d(\xi)$,
see
\cite{DoRo90,ILP,Mu} and \cite[I Sect. 7]{longo 89-90}. Moreover, the 
vacuum representation of $\A$ is associated to the trivial one dimensional 
representation of $G$ and the corresponding Hilbert space $\H_\A$ 
coincides with the subspace of $G$-invariant vectors of $\H_\B$. 

We denote by $\Delta\equiv \Delta_\B$ the semigroup of DHR 
endomorphism of $\A_0$ which are unitarily equivalent to a finite direct 
sum 
of representations $\pi_\xi, \xi\in \hat{G}$. Then the (DHR) braiding on 
$\Delta$ gives in fact a permutation symmetry, $\Delta$ is a dual of
$G$ 
in the sense of Doplicher-Roberts duality theory \cite{DoRo89a,DoRo89b} 
and one can recover the local extension $\B$ by Doplicher-Roberts 
reconstruction theorem \cite{DoRo90}, see \cite[Prop. 3.8]{Mu} for the 
necessary adaptations to conformal nets on $\s1$. 

More generally, let $\A$ be a conformal net on $\s1$ 
and let $\Delta$ be a semigroup of DHR endomorphisms of 
$\A_0$, all covariant with finite dimension.
Assume that the DHR braiding on $\Delta$ is in fact a 
permutation symmetry (para-Bose statistics for the endomorphisms in 
$\Delta$) and that $\Delta$ is specially directed in the sense
of \cite[Sect. 5]{DoRo89a}. Then Doplicher-Roberts construction provides 
a local extension $\B$ of $\A$ and a strongly compact group $G$ of vacuum 
preserving internal symmetries of $\B$ such that $\A$ coincides
with the fixed point net $\B^G$. Note that by \cite[Theorem 3.4]{DoRo89b} 
(see also \cite[Lemma 3.7]{DoRo90})
if $\Delta$ has direct sums, subobjects and conjugates then it is 
specially directed.

In the following we shall use the notation 
$\B = \A \rtimes \Delta$ for the net obtained through the above 
Doplicher-Roberts cross product construction. 

The decomposition in Eq. (\ref{reduction1}) suggests the following 
generalization of the local extensions with compact group action discussed
above, cf. \cite[Sect. 5]{LR}. 

\begin{definition}
\label{defcompacttype}
 A local extension $\B$ of a conformal net $\A$ is of {\bf compact 
type} if the corresponding representation $\pi$ of $\A$ on $\H_\B$ 
satisfies

\begin{equation}
\label{reduction2}
\pi = \bigoplus n_i\pi_i,
\end{equation}
where the $\pi_i$ are (necessarily covariant with positive energy)
mutually inequivalent irreducible subrepresentations 
of $\pi$ with finite dimension appearing  with multiplicity $n_i$ and
having 
finite statistical dimension. 
\end{definition}

Although we did not assume in  Definition \ref{defcompacttype} 
any bound on the  multiplicities $n_i$,  
these turn out to be finite as a consequence of the following proposition, 
cf. \cite[Prop.2.3]{KL} and \cite{CC3} for related results. 
\begin{proposition}
\label{irrlocext}
Let $\B$ be a local extension of compact type of a conformal net $\A$ on 
$\s1$ and let $\pi$ be the corresponding representation of $\A$ on
$\H_\B$.
Then the following hold:
\begin{itemize}
\item[(a)] On $\H_\B$ we have 
\begin{equation} 
\A(I) \vee \A(I^c) = \A(S^1), \quad I\in\I.
\end{equation}
  
\item[(b)] The local extension $\B$ is irreducible, namely 
\begin{equation}
\A(I)' \cap \B(I) = \CC 1, \quad I\in\I. 
\end{equation} 

\item[(c)] Every irreducible representation of $\A$ is contained in $\pi$ 
with finite (possibly zero) multiplicity.
\end{itemize}
\end{proposition}

\begin{proof}
 Let $\theta$ be a representation of $\A$ localized
in $I\in \I$ and equivalent to $\pi$. Then for $J\supset I$, 
$\theta_J$ is a dual canonical endomorphism for the subfactor 
$\A(J)\subset \B(J)$. By  assumption $\pi$ is a direct sum of covariant 
representations with finite statistical dimension. Hence we can find 
isometries $V_i \in \A(I),$ $i\in \NN$, with orthogonal ranges, satisfying
$E_i := V_i{V_i}^* \in \theta(\A(\s1))'$, $\sum_{i\in \NN} E_i=1$
and such that the  representations $\sigma^i$ defined by 
${\sigma^i}_J(\cdot) ={V_i }^*\theta_J (\cdot) V_i$, $J\in \I$ 
are irreducible, covariant, localized in $\I$ and with finite statistical 
dimension. If $T \in \theta_I (\A(I))' \cap \A(I)$  then 
${V_i}^*TV_j{\sigma^j}_I(\cdot)={\sigma^i}_I(\cdot){V_i}^*TV_j$ for 
$i, j \in \NN$ and hence by the equivalence of local and global
intertwiners
for localized representations with finite dimension 
\cite[Theorem 2.3]{GuLo96} we have 
${V_i}^*TV_j{\sigma^j}_J(\cdot)={\sigma^i}_J(\cdot){V_i}^*TV_j$ for every
$J\in \I$. It follows that $E_iT E_j \in  \theta(\A(\s1))'$ and  hence 
$T \in  \theta(\A(\s1))'$.
Since $T\in \theta_I (\A(I))' \cap \A(I)$ was arbitrary and 
$\theta(\A(\s1))' \subset \theta_I (\A(I))' \cap \A(I)$ by Haag duality, 
we conclude that  
$\theta(\A(\s1)) =\theta_I (\A(I)) \vee \theta_{I^c}(A(I^c))$. Hence 
$\pi(\A(\s1)) =\pi_I (\A(I)) \vee \pi_{I^c}(A(I^c))$ which proves (a).
Now, recalling that $U(\psl2) \subset \A(\s1)$, by definition 
of local extensions, we find $\CC 1 = \A(\s1)' \cap \B(I)$ and hence 
(b) follows from (a) and locality.  Finally (c) follows from 
\cite[p. 39]{ILP}
\end{proof}

Since the defining extensions of the fixed point nets under compact groups 
of internal symmetries and the finite index extensions are of compact type 
we can conclude that  (b) of Prop. \ref{irrlocext} generalizes the 
irreducibility results for conformal subsystems in 
\cite[Prop. 2.1]{Car99a} and \cite[Corollary 2.7]{D'ALR} 
(the latter in the local case).

\begin{remark} If $\B$ is a local extension of compact type of a conformal 
net $\A$ on $\s1$ then it follows from Proposition \ref{irrlocext} 
(and its proof) that $\A(I) \subset \B(I),$ $I\in \I$ is an irreducible 
discrete inclusion of infinite factors in the sense of 
\cite[Definition 3.7]{ILP}. 
\end{remark}

We now consider the Virasoro net $\A_{(\rm Vir, 1)}$ with $c=1$. 
By \cite[Prop. 4]{Reh} $\A_{(\rm Vir, 1)}$ is the fixed point net under 
the action of $\so3$ on the conformal net $\A_{\su2_1}$ associated to the 
level one vacuum representation of the loop group $L\su2$. 
The corresponding 
representation $\pi$ of $\A_{(\rm Vir, 1)}$ on $\H_{\A_{\su2_1}}$ 
satisfies 
\begin{equation}
\label{reductionsu21}
\pi = \bigoplus_{j \in \NN_0} (2j+1)\pi^{1}_{j^2}, 
\end{equation}
where $\pi^{1}_{j^2}$ is the representation of $\A_{(\rm Vir, 1)}$ with 
lowest weight $j^2$.  As a consequence $d(\pi^{1}_{j^2})=2j+1$ for each
$j\in \NN$ \cite[Corollary 6]{Reh}. 

We can consider the permutation symmetric semigroup $\Delta$ of DHR 
endomorphisms of $\A_{(\rm Vir, 1)}$  which are 
localized in some $I \in \I_0$ and equivalent to a finite direct sum of  
representations of the type $\pi^{1}_{j^2}$, $j\in \NN_0$. Then, as 
discussed above, $\A_{\su2_1}$ can be identified with 
the Doplicher-Roberts cross product $\A_{(\rm Vir, 1)} \rtimes \Delta$. 

Now let $\B$ be a local extension of compact type of 
$\A_{(\rm Vir, 1)}$ and let $\pi$ be the corresponding representation of 
$\A_{(\rm Vir, 1)}$ on $\H_\B$.   By Prop. \ref{allvirrep} and  
\cite[Theorem 4.4]{Car02} 
every irreducible subrepresentation of $\pi$ is equivalent to 
a DHR endomorphism in $\Delta$ (note that only subrepresentations 
with integer lowest conformal energy can appear) and hence 
\begin{equation}
\label{reductionB}
\pi \simeq \bigoplus_{i\in \NN} \sigma^i,
\end{equation}
where $\sigma_i \in \Delta$, 
for each $i \in \NN$. The local extensions of $\A_{(\rm Vir, 1)}$
with the above property have been independently classified 
by the author (cf. the announcement in \cite{K}) and by Feng Xu
\cite[Sect. 4.2.2]{xu 2003}. The resulting possibilities
are described in the following theorem (we outline our original proof below).

\begin{theorem}
\label{classc=1}
A local extension $\B$ of $\A_{(\rm Vir, 1)}$ is of 
compact type if and only if $\B$ is isomorphic to a fixed point
net ${\A^H}_{\su2_1}$ for some closed subgroup $H$ of $\so3$. 
\end{theorem}
\begin{proof}
The ``if part" is a straightforward consequence of 
the fact that $\A_{\su2_1}$ is an extension of compact type of  
$\A_{(\rm Vir, 1)}$, cf.  Eq. (\ref{reductionsu21}). Now let $\B$ 
be an extension of compact type of $\A_{(\rm Vir, 1)}$ and let $\pi$ 
be the corresponding representation of $\A_{(\rm Vir, 1)}$ on $\H_\B$. 
Given a representation $\theta$ of $\A_{(\rm Vir, 1)}$ 
localized in $I\in\I_0$ and unitarily equivalent to $\pi$ (so that if 
$J \supset I$, $\theta_J$ is a dual canonical endomorphism for  
the inclusion $\A_{(\rm Vir, 1)}(J) \subset \B(J)$) we deduce from 
Eq. (\ref{reductionB}) that $\theta$ is equivalent to (possibly infinite) 
direct sum of DHR endomorphisms in the permutation symmetric 
semigroup $\Delta$ defined after Eq. (\ref{reductionsu21}). It follows 
that the monodromy operator 
$\epsilon_M(\rho,\theta):=\epsilon(\rho, \theta) \epsilon(\theta, \rho)$ 
is trivial (i.e. equal to 1) for every $\rho \in \Delta$. We now use the 
extension of DHR endomorphisms as defined in \cite[Prop. 3.9]{LR} 
(cf. also \cite[Sect. 3.4.7]{Roberts89}) and called $\alpha$-induction 
in \cite{BE}. For every $\rho \in \Delta$, the triviality of the monodromy 
operator $\epsilon_M(\rho,\theta)$ implies that its extension 
$\alpha_\rho$ (we use the notation in \cite{BE}) to $\B_0$ is still 
localized in an interval in $\I_0$, see \cite[Prop. 3.9]{LR}. 

Now the crucial point is that the functorial properties $\alpha$-induction  
(called homomorphism properties in \cite{BE}) imply that 
$\alpha_\Delta$ is still a specially directed permutation symmetric 
semigroup of M\"obius covariant (bosonic) DHR endomorphisms of $\B_0$. 
These functorial properties have been established in \cite{Conti,CDR} for 
inclusions of local nets on the four dimensional Minkowski space-time
and in \cite{BE} for finite index nets of subfactor on the real line. 

Due to the triviality of the monodromy (which is automatic in four 
space-time dimensions) one can use the arguments in \cite{Conti,CDR}
(see also \cite[Sect. 2]{CC2} for an overview) to get the desired
structure 
on $\alpha_\Delta$. Hence, as recalled at the beginning of this section,
we can use the Doplicher-Roberts cross product construction to define 
a local extension $\B \rtimes \alpha_\Delta$ of the conformal net $\B$.

The next point is that the proof of \cite[Theorem 3.5]{CDR} can be adapted 
to our situation to show that there is a natural inclusion 
(up to isomorphism)
$$\A_{\su2_1}=\A_{(\rm Vir, 1)} \rtimes \Delta 
\subset \B \rtimes \alpha_\Delta$$ 
(compatible with $\A_{(\rm Vir, 1)} \subset \B$ ) and in fact it turns out
that $\B \rtimes \alpha_\Delta$ is a local extension of $\A_{\su2_1}$. 
But the latter conformal net has no proper local extensions 
(see e.g. \cite{BMT} and \cite{KL}) and hence we conclude that 
$\B \rtimes \alpha_\Delta =\A_{\su2_1}$. Accordingly $\B = \A^H_{\su2_1}$ 
for some closed subgroup $H$ of $\so3$ as claimed. 
\end{proof}
The above proof relies on specific properties of the 
representation category net $\A_{(\rm Vir, 1)}$, namely on the fact that 
the subcategory of representations with finite statistical dimension
(and trivial univalence) is permutation symmetric, a fact 
that appears to be rather exceptional for conformal nets on $\s1$. 
However in the case of local nets on the four dimensional space-time 
similar ideas have been used by the author and R. Conti 
to study local extensions in a fairly general context \cite{CC3}. As  
matter of fact the above mentioned investigation in \cite{CC3} inspired 
our proof of Theorem \ref{classc=1}. 

Coming back conformal nets on $\s1$ we remark that  there are  well known 
local extensions of the Virasoro net with $c=1$ which are not conformal 
subsystems of $\A_{\su2_1}$ (see e.g. \cite{BMT} ) and hence are not of
compact type as a consequence of Theorem \ref{classc=1}. 
However F. Xu has made further progress and classified the local extensions 
$\B$ of the $c=1$ Virasoro net such that the corresponding 
representation of $\A_{(\rm Vir, 1)}$ on $\H_\B$ contains a subrepresentation 
equivalent to some $\pi^{1}_{j^2}$, $j\in \NN$
\cite[Theorem 4.6]{xu 2003}. The above condition is called 
``spectrum condition" in \cite{xu 2003} where it is conjectured that all 
nontrivial extensions of the $c=1$ Virasoro net have to satisfy it. This 
motivates the following definition
\begin{definition}
A local extension $\B$ of a conformal net $\A$ is {\bf maximally 
non-compact} if the corresponding representation $\pi$ of $\A$ on 
$\H_\B$ satisfies the following condition: the only subrepresentation of 
finite statistical dimension of $\pi$ is the vacuum subrepresentation. 
\end{definition}

From the previous discussion we can conclude that a local 
extension of the $c=1$ Virasoro net satisfies Xu's spectrum condition if and 
only if it is not maximally non-compact. No examples of maximally 
non-compact extensions of this net seem to be known. 
We shall  however exhibit in Section \ref{secinfdim} various examples of 
maximally non-compact extensions for the Virasoro nets with $c>1$.
\medskip

We conclude this section with the following proposition. 

\begin{proposition}
\label{extvir}
Let $\B$ be a local extension of the Virasoro net $\A_{({\rm Vir},c)}$.
Then the following hold:
\begin{itemize}
\item[(a)] $\A_{({\rm Vir},c)}(I)'\cap \B(I)=\CC1$, for every $I\in\I$;
\item[(b)] The net $\B$ is diffeomorphism covariant.
\end{itemize}
\end{proposition}

\begin{proof} Let $\pi$ be the representation of $\A_{({\rm Vir},c)}$ on 
$\H_\B$ associated with the local extension $\B$. If $V$ is the 
corresponding strongly continuous projective unitary representation
of $\diff $ on $\H_\B$ given by (b) of Prop.\ref{directsum}, then  
$V(\gamma) \in \A_{({\rm Vir},c)}(I)$ if $\gamma \in {\rm Diff}(I)$, 
for each $I\in \I$. Moreover, for every $g\in \psl2$ we have  
$V(g)=p(U(g))$, where $U$ is the representation makes $\B$ M\"obius 
covariant. 
Hence, it follows from 
\cite[Theorem 12]{koester03b} that 
$$ \A_{({\rm Vir},c)}(I)'\cap \B(I)= U(\psl2)'\cap \B(I)= \CC1$$
which proves (a).

Now let $I$ be a given interval in $\I$ and let $\gamma \in \diff $ be 
such that $\gamma I = I$. Since $\gamma$ preserves the orientation it must keep
fixed the boundary points of I. An elementary argument (which we omit 
here) then shows that for 
every  $J\in \I$ containing the closure of $I$ we 
can find a diffeomorphism $\gamma^J \in {\rm Diff}(J)$ with 
$\gamma^J |_I = \gamma|_I$, i.e. $\gamma^{-1} \gamma^{J}\in {\rm 
Diff}(I^c)$. Since 
$V(\gamma^{-1} \gamma^{J})\in \A_{({\rm Vir},c)}(I^c)
\subset \B(I^c),$
we find  
$$V(\gamma^J)\B(I) V(\gamma^J)^*=V(\gamma)V(\gamma^{-1} \gamma^{J})
\B(I) V(\gamma^{-1} \gamma^{J})^*V(\gamma)^*=
V(\gamma)\B(I) V(\gamma)^*$$
and hence $V(\gamma)\B(I) V(\gamma)^*\subset \B(J),$ for every $J\in \I$ 
containing the closure of $I$. Thus, being conformal nets continuous from 
outside, we conclude that $V(\gamma)\B(I) V(\gamma)^* \subset \B(I).$  
If $\gamma$ is arbitrary we can always find a  
$g\in \psl2$ such that $gI=\gamma I$. It follows that 
\begin{eqnarray*}
V(\gamma)\B(I) V(\gamma)^* & = & U(g)V(g^{-1} \gamma)\B(I) 
V(g^{-1}\gamma)^* U(g)^* \\
& \subset & \B(gI) = \B(\gamma I),  
\end{eqnarray*}
and hence, for every $I\in \I$, $\gamma \in \diff $, we have 
$$ V(\gamma)\B(I)V(\gamma)^* = \B(\gamma I)$$
and also (b) is proved.

\end{proof}   

\begin{remark} If $\B$ is a diffeomorphism covariant net on $\s1$ and 
$V$ is the corresponding projective unitary representation of $\diff $ one 
can define a covariant subsystem $\C$ of $\B$ by 
\begin{equation}
\C(I)= \{V(\gamma): \gamma \in {\rm Diff}(I)\}'' \quad I\in \I.
\end{equation}
Arguing as in the proof of Prop. \ref{allvirrep} it can be shown that 
the conformal net $\C$ on $\s1$ is isomorphic to $\A_{({\rm Vir},c)}$ for 
some $c\in D\cup [1,\infty)$. It follows that the correspondence between 
diffeomorphism covariant nets on $\s1$ and local extensions of the 
Virasoro nets is one-to-one, cf. \cite{KL}.
\end{remark} 

\section{On the oscillator representations of the Virasoro nets with
$c>1$} 
\label{oscillator}

Let $(\A_\u1, U,\H_\u1)$ be the conformal net generated by 
the $\u1$ chiral current algebra, see \cite{BMT,BS-M}. The Hilbert space 
$\H_\u1$ and the net $\A_\u1$ can be identified with the Fock space 
$e^{\H_1}$, 
where $\H_1$ is acted on by the irreducible representation of $\psl2$ of 
lowest weight 1, and with the corresponding second quantization 
net respectively \cite{GLW}. 

We denote  $\H^{fin}_\u1$ the dense subspace 
of finite energy vectors, i.e. the algebraic direct sum of the $L_0$ 
eigenspaces. Then $\H^{fin}_\u1$ carries the unique irreducible 
lowest weight representation of the oscillator (Heisenberg) algebra 
\begin{eqnarray}
\label{lieu1}
[J_n,J_m]=n\delta_{n+m,0} \quad m, n\in \ZZ, \\
{J_0}=q1,
\end{eqnarray} 
with lowest weight $q=0$, see \cite{BMT} and \cite[Sect. 2.2]{KaRa}. 
The corresponding 
lowest weight vector is the vacuum vector $\Omega$ and for 
$\xi, \psi \in \H^{fin}_\u1, n\in \ZZ$ we have 
\begin{equation}
(\xi, J_n \psi)= (J_{-n}\xi, \psi), 
\end{equation}
(hermiticity).
Note that defining $J^q_n:=J_n, J^q_0=q1$ we obtain a unitary 
representation of the oscillator algebra with arbitrary lowest weight 
$q\in \RR$.

The $\u1$ current $J(z), z=e^{i\vartheta}\in S^1$ is defined as an 
operator valued  distribution  by 
\begin{equation}
J(z)=\sum_{n\in \ZZ} J_n z^{-n-1}
\end{equation}    
and the common invariant domain for the smeared field operators 
$$ J(u) = \oint_\s1\frac{dz}{2\pi i}J(z)u(z)\;\; u\in C^\infty(\s1)$$ 
can be chosen to be the subspace $C^\infty(L_0)$ of smooth vectors 
for $L_0$. For a real function $u\in C^\infty(\s1)$, $J(u)$ is essentially 
self-adjoint and the unitary operators $W(u):=e^{iJ(u)}$ with  
$u\in C^\infty(\s1)\; {\rm  real},\; {\rm supp}\; u \subset I$ generate
$\A_\u1 (I)$ for every $I\in \I$. Moreover the Weyl relations hold:
\begin{equation}
W(u)W(v)=W(u+v){\rm e}^{-\oint_\s1\frac{dz}{4\pi i}u'(z)v(z)}, 
\end{equation}
for real smooth functions $u,v$, where $u'(z)$ denotes the derivative 
$\frac{d}{dz}u(z)=-ie^{-i\vartheta}\frac{d}{d\vartheta}u(e^{i\vartheta})$.

As shown in \cite{BMT} (see also \cite{BMT2}) for every 
$q\in \RR$ there is a covariant irreducible representation 
of $\A_\u1$
(BMT-automorphism) $\gamma_q$ on $\H_\u1$ such that 
\begin{equation}
\label{gammaq}
{\gamma_q}_I(W(u))=e^{ iq\oint_\s1\frac{dz}
{2\pi i}z^{-1}u(z)}W(u)=e^{iJ^q(u)},
\end{equation}
for $I\in \I$, $u\in C^\infty(\s1)$ with support in $I$.  Here the field 
$J^q (z)$ is defined by 
\begin{equation}
J^q(z)=\sum_{n\in \ZZ} J^q_n z^{-n-1}=J(z)+qz^{-1}.
\end{equation}    
$\gamma_{q_1}$ and $\gamma_{q_2}$ are inequivalent if $q_1\neq q_2$. 
Moreover, if $\varphi$ is a real smooth function such that 
$-i\varphi'(z)=z^{-1}q$ for $z\in I$ then 
\begin{equation}
\label{defgamma}
{\gamma_q}_I(\cdot)= {\rm Ad}W(-\varphi) (\cdot),
\end{equation}
and hence $\gamma_q$ is locally implementable by Weyl unitaries. 
In fact Eq. (\ref{defgamma}) can be used to define the representation 
$\gamma_q$. 
Note that $\gamma_0$ is the vacuum representation of $\A_\u1$ 
and that 
for every $I\in\I$, we have 
\begin{equation}
\label{autgamma}
{\gamma_q}_I(\A_\u1(I))=\A_\u1(I). 
\end{equation}

We now come to the oscillator representations of the Virasoro algebra.
For $\lambda, q \in \RR, n\in \ZZ$ the operators 
\begin{equation}
\label{Llq}
L^{(\lambda,q)}_n =
\delta_{n,0}\frac{\lambda^2}{2} +
\frac{1}{2}\sum_{j\in \ZZ}:J^q_{-j}J^q_{j+n}: +i\lambda nJ^q_n,
\end{equation}
where the colons denote normal ordering, 
define a positive energy unitary representation $R(\lambda,q)$
of the Virasoro algebra on  $\H^{fin}_\u1$ with 
central charge $c=1+12\lambda^2$ , see e.g. \cite[Sect. 3.4]{KaRa}. 
Since $L_0$ coincides with $L_0^{(0,0)}$ (by the Sugawara formula)
we have $L^{(\lambda,q)}_0=L_0 + (\lambda^2 + q^2)/2$ and hence 
$\Omega$ is a lowest energy vector for these representations with energy
$(\lambda^2 + q^2)/2$.  

We associate to the above representations
the energy-momentum tensors $T^{(\lambda,q)}(z)$ defined by
\begin{equation}
\label{Tlq0}
T^{(\lambda,q)}(z)=\sum_{n\in \ZZ} L^{(\lambda,q)}_n z^{-n-2}.
\end{equation}
Then the following holds (see \cite[Remark 4.2]{FST})
\begin{equation}
\label{Tlq1}
 T^{(\lambda,q)}(z)=\frac{1}{2}:J^q(z)^2:
 -i\lambda \left(\frac{1}{z} +\frac{{\rm d}}{{\rm d}z} \right)J^q(z) 
 +\frac{\lambda^2}{2z^2}, 
 \end{equation}
 and hence, recalling that $J^q(z)= J(z) +qz^{-1}$, 
 \begin{equation}
 \label{Tlq2}
 T^{(\lambda,q)}(z)=\frac{1}{2}:J(z)^2: +\frac{q}{z}J(z)
 -i\lambda \left(\frac{1}{z} +\frac{{\rm d}}{{\rm d}z} \right)J(z) 
 +\frac{\lambda^2+q^2}{2z^2}.
\end{equation}

For $f\in C^\infty(\s1)$ the smeared field operator 
\begin{equation}
T^{(\lambda,q)}(f)= \oint_\s1\frac{dz}{2\pi i}T^{(\lambda,q)}(z)f(z)
\end{equation}
is well defined on the domain $C^\infty(L_0)$ and leave it globally 
invariant. 
Moreover we see from Eq. (\ref{Tlq2}) that the field 
$T^{(\lambda,q)}(z)$ is local with respect to $J(z)$ in the sense that 
if $f,u \in C^\infty (S^1)$ have disjoint supports the operators 
$T^{(\lambda,q)}(f)$ and $J(u)$ commute on $C^\infty(L_0)$ and that 
$T^{(\lambda,q)}(f)$ is hermitian if $f$ is a real function. 
Finally, it follows  from \cite[Sect. 2]{BS-M} (cf. also \cite{GoWa}) 
that $T^{(\lambda,q)}(f)$ is essentially self-adjoint for each real valued 
smooth function $f$  and that in this case $ e^{i T^{(\lambda,q)}(f)}$ 
commutes with $W(u)$ if the support of the real function $u$ is disjoint from 
the one of $f$.

We now define an isotonous net $\B^{(\lambda,q)}$ on $\H_\u1$ by
\begin{equation}
\label{Blq}
\B^{(\lambda,q)} (I)= 
\{e^{iT^{(\lambda,q)}(f)}: f\in C^\infty(\s1),\; 
{\rm real},\; {\rm supp}\; f\subset I\}'',
\end{equation}
for $I\in \I$.
As a consequence of the above discussion and of Haag duality for $\A_\u1$ 
we obtain the following proposition. 
\begin{proposition}
\label{propBlq}
For every $I\in \I$ we have  
\begin{equation}
\B^{(\lambda,q)} (I) \subset \A_\u1 (I). 
\end{equation}
\end{proposition}
The net $\B^{(\lambda,q)}$ so defined it is not  in general a conformal 
subsystem of $\A_\u1$. And in fact it can be shown that 
$\B^{(\lambda,q)}$
transforms covariantly  with respect to the representation $U$ making 
$\A_\u1$ M\"obius covariant only for 
$(\lambda,q)=(0,0)$, $\B^{(0,0)}$ being the $(c=1)$ Virasoro subnet of 
$\A_\u1$. However, the equality 
$L^{(\lambda,q)}_0=L_0 +(\lambda^2 + q^2)/2$ implies rotation 
covariance for every $(\lambda, q)\in \RR^2$, namely
\begin{equation}
U(r(\vartheta))\B^{(\lambda,q)} (I) U(r(-\vartheta))=
\B^{(\lambda,q)} (r(\vartheta)I)
\quad \vartheta \in \RR, I\in \I. 
\end{equation}

We shall need the following two lemmata.

\begin{lemma}
\label{gammalemma}
For every pair $(\lambda, q)\in \RR^2$ and every $I \in \I$ the following 
holds:
\begin{equation}
 \B^{(\lambda,q)}(I)= {\gamma_q}_I ( \B^{(\lambda,0)}(I)).
\end{equation}
\end{lemma}
\begin{proof}
Let $\varphi, f\in C^{\infty}(S^1)$ be real functions such that 
$-i\varphi'(z)=z^{-1}q$ for $z\in I$ and ${\rm supp} f \subset I$. 
For $\xi, \psi \in \H^{fin}_\u1$ a straightforward calculation shows that 
$$(\xi, W(-\varphi) T^{(\lambda,0)}(f) \psi) =
(T^{(\lambda,q)}(f)\xi, W(-\varphi) \psi).$$

Since $\H^{fin}$ is a common core for $T^{(\lambda,0)}(f)$
and $T^{(\lambda,q)}(f)$ it follows that 
$$W(-\varphi) e^{i T^{(\lambda,0)}(f)}W(\varphi)= 
e^{iT^{(\lambda,q)}(f)}$$
and hence, recalling Eq. (\ref{defgamma})
$${\gamma_q}_I(e^{i T^{(\lambda,0)}(f)})= e^{i T^{(\lambda,q)}(f)}, $$ 
cf. \cite[p. 123]{BS-M} and \cite[p. 361]{BMT2}. The conclusion 
then follows from the definition of $\B^{(\lambda,q)} (I)$ given in Eq. 
(\ref{Blq}). 
\end{proof}

\begin{lemma} 
\label{Rlq}
The representation $R(\lambda,q)$ is irreducible for every 
$\lambda\neq 0$ and $q\in \RR$. 
\end{lemma}
\begin{proof}
The character $\chi_{(\lambda,q)}(t)$, 
$t\in (0,1)$ of the representation $R(\lambda,q)$ is given by 
$$\chi_{(\lambda,q)}(t)={\rm Tr}( t^{L^{(\lambda,q)}_0 })= 
t^{\frac{\lambda^2 + q^2}{2} } p(t),$$
where $p(t)=\prod_{n=1}^{\infty} (1-t^n)^{-1}= 
{\rm Tr}( t^{L_0 })$ and hence the conclusion follows 
since, by \cite[Eq. (3.15) and Prop. 8.2]{KaRa}, it coincides with 
the character of the irreducible representation $L(c,h)$ of Vir 
with central charge $c=1+12\lambda^2$ and lowest weight 
$h=(\lambda^2 + q^2)/2$. 
\end{proof}

\begin{corollary}
\label{corollaryh}
Let $\A_{({\rm Vir},c)}$ be the Virasoro net with 
central charge $c=1+12\lambda^2, \lambda\neq 0$ and let $\pi^{c}_{h}$ 
be the (irreducible) representation of $\A_{({\rm Vir},c)}$ with lowest 
weight $h=(\lambda^2 + q^2)/2$ as defined in Subsect. \ref{virsubsec}. 
Then there is a representation $\pi_{(\lambda,q)}$ of 
$\A_{({\rm Vir},c)}$ on 
$\H_\u1$, unitarily equivalent to $\pi^{c}_{h}$, such that for every 
$I\in \I$ the following holds:
\begin{equation}
{\pi_{(\lambda,q)}}_I (\A_{({\rm Vir},c)}(I) )=\B_{(\lambda,q)}(I). 
\end{equation}
\end{corollary} 
We shall need the following proposition in the next section.  
\begin{proposition}
\label{constdim} 
Let $\A_{({\rm Vir},c)}$ be a Virasoro net with $c>1$. Then, 
if $h \geq (c-1)/24$ we have $d(\pi^{c}_{h})=d(c)$, where 
$d(c) \in [1, \infty]$ does not depend on $h$ and satisfies 
$d(c)>1$. 
\end{proposition}
\begin{proof} 
The assumption on the range of $c$ and $h$ implies that we can find 
$\lambda \neq 0$ and $q\in \RR$ such that $c=1+12\lambda^2$ and 
$h=(\lambda^2 + q^2)/2$. Then it follows from 
Corollary \ref{corollaryh} that $d(\pi^{c}_{h})=d(\pi_{(\lambda,q)})$ and 
we have to show that the latter does not depend on $q$. 

By  Eq. (\ref{def.dim.}) and
Corollary \ref{corollaryh} we find 
$$d(\pi_{(\lambda,q)})^2=[\B_{(\lambda,q)}(I^c)': \B_{(\lambda,q)}(I)],
\quad 
I\in \I. $$
From Proposition \ref{propBlq} and Haag duality for $\A_\u1$ it 
follows that
$$\B_{(\lambda,q)}(I) \subset \A_\u1 (I)=\A_\u1 (I^c)' \subset 
\B_{(\lambda,q)}(I^c)'$$ 
and hence, using the multiplicativity of the minimal index 
\cite{longo 92} (cf. the proof of \cite[Prop. 3.1]{Car02}), 
that $$d(\pi_{(\lambda,q)})^2= [\A_\u1(I): \B_{(\lambda,q)}(I)]
\cdot [\A_\u1(I^c): \B_{(\lambda,q)}(^c)].$$
Now, using Lemma \ref{gammalemma} and Eq. (\ref{autgamma}) we
find, for an arbitrary  $J\in I$, 
\begin{eqnarray*}
[\A_\u1(J): \B_{(\lambda,q)}(J)] & = & 
[{\gamma_q}_J (\A_\u1(J) ): {\gamma_q}_J (\B_{(\lambda,0)}(J)) ]
\\
& = & [\A_\u1(J): \B_{(\lambda,0)}(J)],
\end{eqnarray*}
and hence 
$$d(\pi_{(\lambda,q)})^2= [\A_\u1(I): \B_{(\lambda,0)}(I)]
\cdot [\A_\u1(I^c): \B_{(\lambda,0)}(I^c)]$$
does not depend on $q$. 

Finally if $d(c)=1$ then, for every $I\in \I$, 
$\pi_{(\lambda,q)}(\A_{({\rm Vir},c)}(I))= \A_\u1(I)$
which is impossible since $\A_\u1$ is strongly additive 
(see \cite{BS-M,GLW}) while $\A_{({\rm Vir},c)}$ it is not.
\end{proof}

\section{Sectors with infinite dimension and maximally 
non-compact local extensions.}
\label{secinfdim}

Let $D\subset [\frac{1}{2}, 1)$ be the set of
allowed values of the central charge in the discrete series 
representations of Vir as defined in 
Subsect. \ref{virsubsec} and   
let $c \in (D+1) \cup [2, \infty)$. Then $c-1$ is an allowed value 
of the central charge and the tensor product net 
$\A_{({\rm Vir},c-1)} \otimes \A_{\su2_1}$  
is a local extension of 
$\A_{({\rm Vir},c)}$. The representation of $\A_{({\rm Vir},1)}$ on 
$\H_{\A_{\su2_1}}$ contains the irreducible lowest weight representation 
$\pi^{1}_{j^2}$, $j\in \NN_0$ with multiplicity $2j+1$ 
(see Eq. (\ref{reductionsu21})) and hence 
the multiplicity $m(c,j)$ of  $\pi^{c}_{j^2}$ in the representation of 
$\A_{({\rm Vir},c)}$ on
$\H_{\A_{({\rm Vir},c-1)}} \otimes \H_{\A_{\su2_1}}$ satisfies 
$m(c,j) \geq 2j +1$ for every $j\in \NN_0$. 
We are now ready to prove the following theorem, cf. 
\cite[Theorem 4.4]{Car02} and the guess in \cite[Sect. 2]{Reh}.

\begin{theorem} 
\label{infdim}
If $c\in (D+1) \cup [2, \infty)$ and 
$h \geq (c-1)/24$ then $d(\pi^{c}_{h})=\infty$.
\end{theorem}

\begin{proof}
Let $\pi$ be the representation of $\A_{({\rm Vir},c)}$ in 
$\H_{\A_{({\rm Vir},c-1)}} \otimes \H_{\A_{\su2_1}}$ as described above. 
Then, as explained in Sect. \ref{preliminaries}, $\pi$ is 
unitarily equivalent to a representation $\theta$ on 
$\H_{\A_{({\rm Vir},c)}}$ localized in an interval $I_0 \in \I$ and for
every 
$I \in \I$ with $I_0 \subset I$ $\theta_I$ is a dual canonical
endomorphism 
for the  inclusion  
$\A_{({\rm Vir},c)} (I) \subset 
\A_{({\rm Vir},c-1)}(I) \otimes \A_{\su2_1} (I),$ which is irreducible 
because of Prop. \ref{extvir}. Now let 
$\rho^{c}_{j^2}$
be a representation of $\A_{({\rm Vir},c)}$ on $\H_{\A_{({\rm Vir},c)}}$, 
unitarily equivalent to  $\pi^{c}_{j^2}$ and localized in $I_0$ and let 
$I\in \I$ be an interval containing $I_0$. As shown just before the 
statement of this theorem the multiplicity $m(c,j)$ of the representation 
$\rho^{c}_{j^2}$ in $\theta$ satisfies $m(c,j) \geq 2j +1$. Hence (by 
Haag duality) the endomorphism ${\rho^{c}_{j^2}}_I$ is contained in 
$\theta_I$ with multiplicity $n(c,j) \geq 2j +1$ for each $j\in \NN_0$. 
Now, it follows from Prop. \ref{constdim} that  
$d(\rho^{c}_{j^2})=d(c)$, for each $j \geq \sqrt{(c-1)/24}$, 
where $d(c)$ does not depend on $j$. Let us assume that $d(c) < \infty$. 
Then by  \cite[Corollary 2.10]{GuLo96}  ${\rho^{c}_{j^2}}_I$ is an 
irreducible endomorphism of  $\A_{({\rm Vir},c)}(I)$ for every 
$j \geq \sqrt{(c-1)/24}$ and by \cite[p. 39]{ILP}
we conclude that $2j+1\leq n(c,j) \leq d(c)^2$ for every 
$j \geq \sqrt{(c-1)/24}$, in contradiction with 
the assumption $d(c) < \infty$. Hence $d(c) =\infty$ and the conclusion 
follows from Prop. \ref{constdim}.
\end{proof}

\begin{corollary} 
\label{indexcorollary}
If $c\in (D+1) \cup [2, \infty)$ and $\B$ is a local 
extension of compact type of $\A_{({\rm Vir},c)}$ then the index 
$[ \B :\A_{({\rm Vir},c)} ]$ is finite. 
\end{corollary}
\begin{proof}
Let $\pi$ the representation of  $\A_{({\rm Vir},c)}$
on $\H_\B$ defined by the local extension $\B$. Only representations with 
integer lowest weigh can appear in the decomposition of $\pi$. But there 
are only a finite number of positive integers $m$ satisfying 
$m < (c-1)/24$ and hence, by Theorem \ref{infdim}, only a 
finite number of irreducible DHR sectors can appear in the decomposition 
of $\pi$. Now, recalling that the inclusion 
$\A_{({\rm Vir},c)}(I) \subset \B(I), \quad I\in \I,$
is irreducible, the conclusion follows from (the proof of) 
\cite[Prop. 2.3]{KL}. 
\end{proof}

Now let $G$ be a simply connected compact Lie group with simple 
Lie algebra ${\rm Lie}(G)$ and let $k$ be a positive integer. We denote 
by $\A_\gk$ the conformal net associated to the vacuum representation 
of the corresponding Loop group (or affine Lie algebra) at level $k$ 
(see \cite{FrG,Reh,Tol97,Was A}). As it is well known, the Sugawara 
formula (see e.g. \cite[Sect. 15.2]{DMS} and 
\cite[Sect.10.1]{KaRa}), 
implies that the net $\A_\gk$ is a local extension of the Virasoro net 
$\A_{({\rm Vir},c)}$ with central charge 
\begin{equation}
c \equiv c(\gk) = \frac{{\rm dim}(G) k}{k + h^{\vee}},
\end{equation}
where $h^{\vee}$ is the dual Coxeter number of ${\rm Lie}(G)$, 
cf. \cite[Sect. III.7]{FrG} and  \cite[Sect.1]{Reh}. 
The central charge $c(\gk)$ is bounded by
\begin{equation}
r\leq c(\gk) \leq {\rm dim}(G), 
\end{equation}
where $r$ is the rank of ${\rm Lie}(G)$ 
and the lower bound is saturated only for a simply laced Lie algebras 
at level $k=1$. Note that $c(\gk) < 2$ implies that $r=1$ and thus that 
$G= \su2$. In the latter case we have $c({\su2}_k)= 3k/(k+2)$.
If $k\geq 4$ we have $c({\su2}_k) \geq 2$. The remaining 
possibilities are $c({\su2}_1)=1$, $c({\su2}_2)=1+{1}/{2}$ 
and $c({\su2}_3)= 1+{4}/{5}$. We summarize the above discussion in 
the following lemma.

\begin{lemma}
\label{cgklemma}
If $\gk \neq {\su2}_1$ then 
$c(\gk) \in (D+1) \cup [2, \infty)$. 
\end{lemma} 

Recall that there is a strongly continuous representation
of $G$ in the (unitary) group of internal symmetries of $\A_\gk$ leaving
the vacuum invariant. This representation is not in general faithful and 
its kernel coincides with the (finite) center $Z(G)$ of $G$. 
It is known that the fixed point net $\A^G_\gk$ satisfies 
\begin{equation}
\A_{({\rm Vir},c)} \subset \A^G_\gk \subset \A_\gk, \quad c=c(\gk),
\end{equation}
see \cite{Reh}. 
In particular, being $G/Z(G)$ infinite, the index 
$[\A_\gk :\A_{({\rm Vir},c)}]$ is infinite.

\begin{corollary} 
\label{noncompactcorollary}
If $\gk \neq {\su2}_1$ then the local extension 
$\A_\gk$ of $\A_{({\rm Vir},c)}$, $c=c(\gk)$, is not of compact type. 
\end{corollary}
\begin{proof}
Due to Lemma \ref{cgklemma} we can apply 
Corollary \ref{indexcorollary} and the conclusion follows from 
$$[\A_\gk :\A_{({\rm Vir},c)}]=\infty. $$
\end{proof}

The following consequence of Corollary 
\ref{noncompactcorollary} has been pointed out by K.-H. Rehren in
\cite{Reh}
with a different argument based on the comparison of characters. 

It can also be proved using \cite[Theorem 2.4]{xu 2003} and the fact
that 
$\A_{({\rm Vir},c)}$ is not strongly additive when $c>1$.

\begin{corollary} 
\label{fixpoint}
If $\gk \neq {\su2}_1$, then the inclusion
$\A_{({\rm Vir},c)} \subset\A^G_\gk$ is proper. 
\end{corollary}

The next result shows that maximally non-compact local extensions 
naturally appear for the Virasoro nets $c>1$. 
\begin{proposition}
\label{gkmaxnon-comp} 
If $\gk \neq {\su2}_1$ and $c=c(\gk) \leq 25$ then 
$\A_\gk$ is a maximally non-compact local extension of 
$\A_{({\rm Vir},c)}$. 
\end{proposition}

\begin{proof}
The representation $\pi$ of $\A_{({\rm Vir},c)}$ in 
$\H_{\A_\gk}$ can only have irreducible subrepresentations with 
a nonnegative integer lowest weight. Since by assumption 
$(c-1)/24 \leq 1$, it follows from Theorem \ref{infdim} that the
only 
subrepresentation $\pi$ with finite dimension is the vacuum
representation. 
Hence the extension is maximally non-compact.
\end{proof}

For $\suN$ $h^\vee=N$ and hence $c(\suN_k)=k(N^2-1)/(N+k)$ and we 
see that Prop. \ref{gkmaxnon-comp} gives an infinite series of 
maximally non-compact extensions of the $c>1$ Virasoro nets. 
Examples are: $\su2_k$, $k>1$; 
${\rm SU}(3)_k, {\rm SU}(4)_k, {\rm SU}(5)_k$, k arbitrary; 
$\suN_1$, $2\leq N \leq 26$. Actually the same proof of 
Prop. \ref{gkmaxnon-comp}, together with Prop. \ref{allvirrep}, gives the 
following stronger result. 

\begin{theorem}
\label{cmaxnon-comp} If $c\in (1+D) \cup [2,25]$ then every local
extension 
of the Virasoro net $\A_{({\rm Vir},c)}$ is maximally non-compact. 
In particular $\A_{({\rm Vir},c)}$ has no local extensions of compact
type. 
\end{theorem}

\appendix

\section{Appendix}
\label{appendix}
In this appendix we give a differentiability result for the 
representations of $\diff$ which is used in the proof of 
Prop. \ref{allvirrep}. This result as been essentially obtained by T. Loke 
\cite{loke} (cf. also \cite{Was A} for analogous results for loop 
groups) and here we consider the necessary modifications we need in this 
paper. 
We shall closely follow the discussion in \cite[Chap. I]{loke}.
  
An element of the group $\mob$ of M\"obius transformations of $\s1$  
is given by a map 
$z\mapsto \frac{\alpha z + \beta}{\overline{\beta}z+\overline{\alpha}}$, 
where $\alpha, \beta$ are complex numbers satisfying 
$|\alpha|^2 -|\beta|^2 =1$. $\mob$ is a Lie subgroup of $\diff$ isomorphic 
to $\psl2$. The corresponding Lie subalgebra of $\vect$ is spanned by the 
vector fields 
\begin{equation}
\label{gensl2}
x:= -\sin{\vartheta}\frac{d}{d\vartheta},\;
y:=-\cos{\vartheta}\frac{d}{d\vartheta},\; h:= \frac{d}{d\vartheta},
\end{equation} 
whose brackets are given by
\begin{equation}
\label{liesl2}
[h,x]=-y,\; [h,y]=x,\; [x,y]=h.
\end{equation}
More generally, for each $n\in \NN$, the vector fields 
\begin{equation}
x_n:= -\frac{1}{n}\sin{n\vartheta}\frac{d}{d\vartheta},\;
y_n:=-\frac{1}{n}\cos{n\vartheta}\frac{d}{d\vartheta},\; 
h_n:= \frac{1}{n}\frac{d}{d\vartheta},
\end{equation}
span isomorphic Lie subalgebras of $\vect$ each associated to a Lie 
subgroup $\mob_n$ of $\diff$. Clearly $\mob_1 =\mob$ and it is not hard to  
see that, for each $n>1$, $\mob_n$ is isomorphic to an n-fold covering of 
$\psl2 \simeq \mob$ and that the corresponding covering map transforms the 
one-parameter 
group $\exp(th_n)$ into the one-parameter subgroup $r(t)$ of rotations of 
$\psl2$.

Now let $V$ be a strongly continuous projective unitary representation of 
$\diff$ on a separable Hilbert space. 
For every $n\in \NN$, the restriction of $V$ 
to $\mob_n$ lifts to a strongly continuous unitary 
representation $U_n$ of $\widetilde{\psl2}$. Note that 
$\exp(2\pi h_n)^n=1$ and 
hence $U_n(\tilde{r}(2\pi))^n = U_n(\exp(2\pi h_n))^n =\chi_n 1$ for a 
suitable complex number $\chi_n$ of modulus one. In particular 
$U(\tilde{r}(2\pi))$ has finite spectrum for each $n\in \NN$.  

Now let $\frac{1}{n}X_n$, $\frac{1}{n}Y_n$ and $\frac{i}{n}(L_0 +c_n)$, 
$c_n\in\RR$, $c_1=0$, be 
the skew-adjoint generators of the one-parameter groups of unitaries
$U_n(\exp(tx_n))$, $U_n(\exp(ty_n))$, and $U_n(\exp(th_n))$, respectively.
On the dense subspace $\D_n\subset \H$ of $C^\infty$ vectors for the 
representation $U_n$ the above operators define a representation 
of the Lie algebra (\ref{liesl2}) and hence we have on $\D_n$
\begin{equation} 
\label{liesl2b}
[iL_0,X_n]=-nY_n,\; [iL_0,Y_n]=nX_n,\; [X_n,Y_n]=in(L_0+c_n),\;n\in\NN.
\end{equation} 

If $V$ is a positive energy representation,
since the unitary operator 
$e^{i2\pi L_0}$ acts as multiplication by a complex number, the spectrum 
of $L_0$ is pure point and every eigenvalue is of the form $h+n$, where 
$h\geq 0$ is the lowest eigenvalue of $L_0$ and $n$ is a nonnegative integer.

Now let $\H^{fin}$ be the linear span of the eigenvectors of $L_0$. 
Loke as shown in \cite[Sect. I.1]{loke} that if a positive energy 
representation $V$ is such that the eigenspaces of $L_0$ are all 
finite-dimensional, then 
\begin{equation}
\H^{fin}\subset \bigcap_{n\in \NN}\D_n.
\end{equation}

Moreover he proved that the operators $L_0$, $L_n:=iY_n-X_n$ and 
$L_{-n}:=iY_n + X_n$, $n\in \NN$ define a unitary representation of Vir on 
$\H^{fin}$ and that the corresponding energy-momentum tensor 
\begin{equation}  
T(z)=\sum_{n\in \ZZ}L_n z^{-n-2}
\end{equation}
extends to an operator valued distribution on the subspace of smooth $L_0$
vectors such that $T(f)$ is essentially self-adjoint on $\H^{fin}$ for 
each $f\in \vect$ and satisfies 
\begin{equation} 
p(e^{iT(f)})=V(\exp(f)).
\end{equation}
The finite dimensionality of the $L_0$ eigenspace is used in \cite{loke} 
to infer that for each $n\in \NN$, the representation $U_n$ is a direct 
sum of positive energy representations of $\widetilde{\psl2}$ and that 
$\D_n \subset \H^{fin}$. 
However these facts hold for every positive energy representation $V$, as 
a consequence of the proposition below (applied to each representation 
$U_n$) and hence the results of Loke 
described above hold (without any essential modification in the proofs) 
also if the finite dimensionality of the eigenspaces of $L_0$ is not 
assumed. Moreover the representation of Vir on $\H^{fin}$ so obtained can 
be seen to be irreducible (and hence unitarily equivalent to some 
$L(c,h)$, cf. \cite[Remark 3.5]{KaRa}) if and only if the corresponding 
projective representation $V$ of $\diff$ is irreducible, 
cf. Lemma 2.2. in \cite[Sect. I.2]{loke}. 

\begin{proposition} Let $U$ be a strongly continuous unitary 
representation of $\widetilde{\psl2}$ on a separable Hilbert space 
$\H$ and let $L_0$ be the 
self-adjoint generator of the restriction of $U$ to the lifting $\tilde{r}(t)$ 
of the one-parameter rotation subgroup of $\psl2$. Assume that the spectrum 
of $L_0$ is bounded from below and that the one of 
$U(\tilde{r}(2\pi))$ is finite. Then the following hold:
\begin{itemize}
\item[(a)] $U$ is a positive energy 
representation (i.e. $L_0$ has a nonnegative spectrum) and it is 
completely
reducible to a direct sum of irreducible subrepresentations.   
\item[(b)] Every eigenvector of $L_0$ is a smooth vector for the 
representation $U$. 
\end{itemize}
\end{proposition}

\begin{proof} 
If $U$ is assumed to be irreducible then the positivity of the energy 
follows from the bound on the spectrum of $L_0$ as a consequence of the 
classification of the irreducible representations of $\widetilde{\psl2}$ 
\cite{pukanzsky} (cf. \cite[Sect. I.1.3]{loke}) and hence the 
positive energy condition for $U$ follows in general by direct integral 
decomposition. Then (a) follows e.g. from \cite[Lemma 8]{koester02}. 
As a consequence there is an increasing sequence $0=n_1<n_2...$ of 
nonnegative integers (which is possibly finite) and a decomposition 
$$\H = \bigoplus_{k}\H_k$$ 
such that the restriction of $U$ to $\H_k$ is a (possibly infinite) 
multiple of an irreducible  representation with lowest weight $h+n_k$. 
Hence if $\psi$ is an eigenvector of $L_0$ corresponding to the 
eigenvalue $\lambda$ we can write 
$$\psi = \sum_{n_k \leq \lambda -h} (\psi_k,\psi) \psi_k$$
where $\psi_k \in \H_k$ is a normalized eigenvector of $L_0$.   

Since every eigenvector of the generator of rotations in an 
irreducible representation of $\widetilde{\psl2}$ is smooth (see e.g. 
\cite[Sect. I.1]{pukanzsky}), each $\psi_k$ is smooth and 
hence $\psi$ is smooth vector for the representation $U$ so that also (b) 
is proved. 
\end{proof}

We can summarize the discussion in this appendix in the following theorem, 
cf. \cite[Sect. I.2.4]{loke}. 
\begin{theorem} 
\label{diffdiff}
Let $V$ be a strongly continuous positive energy 
projective unitary irreducible representation of $\diff$ on a 
(necessarily separable) Hilbert space $\H$. Then $V$ is unitarily 
equivalent to the unique projective unitary representation $V_{(c,h)}$ 
which integrates the Vir-module $L(c,h)$ for some $c>0, h\geq 0$. In 
particular the corresponding generator of rotations $L_0$ has 
finite-dimensional eigenspaces.
\end{theorem}

\medskip
\noindent{\bf Acknowledgements.} The author would like to thank R. Conti,
S. K\"oster and R. Longo for discussions, explanations and comments. 
Theorem \ref{classc=1} has been announced at the Miniworkshop 
``Conformal Field Theory. An Introduction" held in Rome on March 2003. 
The author thanks the organizers D. Guido and (again) R. Longo for the 
invitation.


\medskip

\begin{thebibliography}{99}

{\scriptsize

\bibitem{BCL} Bertozzini P., Conti R., Longo R.:
Covariant sectors with infinite dimension and positivity of the energy.
{\it Comm. Math. Phys.} {\bf 193} (1998), 471--492.

\bibitem{BE} B\"{o}ckenhauer J., Evans D.E.: Modular invariants graphs
and $\alpha$-induction for nets of subfactors I. {\it Comm. Math. Phys.}
{\bf 197} (1998), 361--386.

\bibitem{BMT} Buchholz D., Mack G., Todorov I.T.: The current algebra on
the
circle as a germ of local field theories. {\it Nucl. Phys. B} (Proc.
Suppl.)
{\bf 5B} (1988), 20--56.

\bibitem{BMT2}  Buchholz D., Mack G., Todorov I.T.: Localized
automorphisms of the U(1)-current. In \cite{Kast},

\bibitem{BS-M} Buchholz D., Schulz-Mirbach H.: Haag duality in conformal
quantum field theory. {\it Rev. Math. Phys.} {\bf 2} (1990), 105--125.

\bibitem{Buch90} Buchholz D.: Introduction to conformal QFT in two 
dimensions. Unpublished manuscript (1990).

\bibitem{Car99a} Carpi S.: Classification of subsystems for the
Haag-Kastler nets generated by $c=1$ chiral current algebras.
{\it Lett. Math. Phys.} {\bf 47} (1999), 353--364.

\bibitem{Car02} Carpi S.: The Virasoro algebra and sectors with infinite 
statistical dimension. \\ math.OA/0203027, to appear in 
{\it Ann. H. Poincar\'e}.  

\bibitem{CC2} Carpi S., Conti R.: Classification of subsystems, 
local symmetry generators and intrinsic definition of local observables.
In: R. Longo (ed.) ,{\it Mathematical physics in mathematics and physics.} 
Fields Institute Communications Vol. 30, AMS, Providence, RI, 2001, 
pp.83--103.

\bibitem{CC3} Carpi S., Conti R. : in preparation.

\bibitem{Conti} Conti R.: Inclusioni di algebre di von Neumann e teoria 
algebrica dei campi. 
{\it Ph.D. Thesis}, Universit\`a di Roma {\it Tor Vergata} (1996).

\bibitem{CDR} Conti R., Doplicher S., Roberts J. E.: 
Superselection theory for subsystems.
{\it Comm. Math. Phys.}, {\bf 218} (2001), 263--281.

\bibitem{D'ALR} D'Antoni C., Longo R., Radulescu F.: Conformal nets, 
maximal temperature and and models from free probability. 
{\it J. Operator Theory} {\bf 45} (2001), 195--208.  
 
\bibitem{DMS} Di Francesco Ph., Mathieu P., S\'en\'echal D.: {\it
 Conformal Field Theory.} Springer-Verlag, Berlin-Heidelberg-New York, 
1996.

\bibitem{DHR0} Doplicher S., Haag, R., Roberts J.E.:
Fields, observables and gauge transformations I, II.
{\it Comm. Math. Phys.} 
{\bf 13} (1969), 1--23;
{\it Comm. Math. Phys.}
{\bf 15} (1969), 173--200.

 \bibitem{DoRo89a} Doplicher S., Roberts J. E.:
 Endomorphisms of $C^*$--algebras, cross products and duality
 for compact groups.
 {\it Ann. Math.} 
 {\bf 130} (1989), 75--119.
 
 \bibitem{DoRo89b} Doplicher S., Roberts J. E.:
 A new duality theory for compact groups.
 {\it Invent. Math.} 
 {\bf 98} (1989), 157--218.

\bibitem{DoRo90} Doplicher S., Roberts J.E.: Why there is a field
algebra with a compact gauge group describing the superselection
structure in particle physics. {\it Comm. Math. Phys.\/}
{\bf 131} (1990), 51--107.

\bibitem{FJ} Fredenhagen K., J\"{o}r{\ss} M.: Conformal
Haag-Kastler nets, pointlike localized fields and the
existence of operator product expansions. {\it Comm. Math. Phys.} 
{\bf 176} (1996), 541--554.

\bibitem{FRS2} Fredenhagen K., Rehren K.-H., Schroer B.: Superselection 
sectors with braid group statistics and exchange algebras II. Geometric 
aspects and conformal covariance. {\it Rev. Math. Phys.} 
{\bf Special Issue} (1992), 113--157.

\bibitem{FST} Furlan P., Sotkov G. M., Todorov I.T.: Two-dimensional 
conformal quantum field theory. {\it Riv. Nuovo Cimento} {\bf 12}, N.6 
(1989) 1--202.

\bibitem{FrG} Gabbiani F., Fr\"{o}hlich J.: Operator algebras and
conformal field theory. {\it Comm. Math. Phys.} {\bf 155} (1993), 
569--640.

\bibitem{GKO} Goddard P., Kent A., Olive D.: Unitary representations of 
the Virasoro and super-Virasoro algebra. {\it Comm. Math. Phys.} {\bf 103} 
(1986), 105--119.

\bibitem{GoWa} Goodman R. and Wallach N. R.: Projective unitary 
positive-energy representations of ${\rm Diff}(\s1)$. 
{\it J. Funct. Anal.} {\bf 63}, (1985) 299--321.

\bibitem{GuLo96} Guido D., Longo R.: The conformal spin and statistic 
theorem. {\it Comm. Math. Phys.} {\bf 181} (1996), 11--35.

\bibitem{GLW} Guido D., Longo R., Wiesbrock H.-W.: Extensions of
conformal
nets and superselection structures. {\it Comm. Math. Phys.} {\bf 192}
(1998), 217--244. 

\bibitem{Haag} Haag R.: {\it Local Quantum Physics.} 2nd ed.\,  
Springer-Verlag,  Berlin-Heidelberg-New York, 1996.

\bibitem{ILP} Izumi M., Longo R., Popa S.: A Galois correspondence 
for compact groups of automorphisms of von Neumann algebras with a
a generalization to Kac algebras. {\it J. Funct. Anal.} {\bf 155} (1998), 
25--63.

\bibitem{jones} Jones V.: Index of subfactors. {\it Invent. Math.} {\bf
72}
(1983), 1--25.

\bibitem{KaRa} Kac V. G., Raina A. K.: {\it Bombay Lectures on Highest 
Weight Representations of Infinite Dimensional Lie Algebras.} World
Scientific, Singapore, 1987.

\bibitem {Kast} 
Kastler D. ed.:
{\it The algebraic theory of superselection sectors.} 
World Scientific, Singapore, 1990.  

\bibitem{K} Kawahigashi Y.: Classification of operator algebraic conformal 
field theories. \\ math.OA/0211141.

\bibitem{KL} Kawahigashi Y., Longo R.: Classification local conformal 
nets. Case $c<1$. math.OA/0211141, to appear in {\it Ann. Math.}

\bibitem{KL2} Kawahigashi Y., Longo R.: Classification of two-dimensional 
local conformal nets with $c<1$ and 2-cohomology vanishing for tensor 
categories. math-ph/0304022.

\bibitem{KLM} Kawahigashi Y., Longo R., M\"{u}ger M.: Multi-interval 
subfactor and modularity of representations in conformal field theory.
{\it Comm. Math. Phys.} {\bf 219} (2001), 631--669.

\bibitem{koester02} K\"oster S.: Conformal transformations as observables.
{\it Lett. Math. Phys.} {\bf 61} (2002), 187--198.
 
\bibitem{koester03a} K\"oster S.: Absence of stress energy tensor in 
${\rm CFT}_2$ models. math-ph/0303053. 

\bibitem{koester03b} K\"oster S.: Local nature of cosets models. 
math-ph/0303054. 

\bibitem{kosaki} Kosaki H.: Extension of Jones' theory on index to
arbitrary subfactors. {\it J. Funct. Anal.} {\bf 66} (1986), 123--140.

\bibitem{loke} Loke T.: Operator algebras and conformal field theory
of the discrete series representation of $\diff$. PhD Thesis, University 
of Cambridge, 1994. 

\bibitem{longo 89-90}  Longo R.: Index of subfactors and statistics of
quantum fields. I. {\it Comm. Math. Phys.} {\bf 126} (1989), 217--247, 
and II. Correspondences, braid group statistics and Jones polynomial. 
 {\it Comm. Math. Phys.} {\bf 130} (1990), 285--309.

\bibitem{longo 92} Longo R.: Minimal index and braided subfactors. 
{\it J. Funct. Anal.} {\bf 109} (1992), 98--112.

 \bibitem{longo 2003} Longo R.: Conformal subnets and
intermediate subfactors. {\it Comm. Math. Phys.} {\bf 237} (2003), 7--30.

\bibitem{LR} Longo R., Rehren K.-H.: Nets of subfactors.
{\it Rev. Math. Phys.} {\bf 7} (1995), 567--597.

\bibitem{mack} Mack G.: Introduction to conformal invariant quantum field 
theory in two and more dimensions. In G. t' Hooft et al. Eds.: 
{\it Non perturbative quantum field theory.} Plenum Press, New York, 1988, 
pp.353--383.

\bibitem{Milnor} Milnor J.: Remarks on infinite-dimensional Lie groups. 
In B.S. De Witt and R. Stora Eds.: {\it Relativity, groups and topology 
II.} Les Houches, Session  XL, 1983, Elsevier, Amsterdam, New York, 1984, 
pp. 1007--1057.

\bibitem{Mu} M\"uger M.: 
On charged fields with group symmetry and degeneracies of
Verlinde's matrix $S$. 
{\it Ann. Inst. H. Poincar\'e}
{\bf 71} (1999), 359--394.

\bibitem{pukanzsky} Puk\'anzsky L.: The Pancherel formula for the 
universal covering group of SL(2,$\RR$). {\it Math. Annalen} {\bf 156}
(1964), 96--143.

\bibitem{Reh} Rehren K.-H.: A new view of the Virasoro algebra.
{\it Lett. Math. Phys.} {\bf 30} (1994), 125--130.

\bibitem{Roberts89} Roberts J.E.:
Lectures on algebraic quantum field theory.
In \cite{Kast}, pp. 1--112.

\bibitem{takesaki} Takesaki M.: {Theory of operator algebras I.} 
Springer-Verlag, Berlin-Heidelberg-New York, 2002.  

\bibitem{Tol97} Toledano Laredo V.: Fusion of positive energy 
representations of $LSpin_{2n}$, PhD Thesis, University of 
Cambridge, 1997. 

\bibitem{Tol99} Toledano Laredo V.:  Integrating unitary representations
of infinite-dimensional Lie groups, 
{\it J. Funct. Anal.} {\bf 161} (1999), 478--508.

\bibitem{Was A} Wassermann A.: Operator algebras and conformal field
theory III: Fusion of positive energy representations of SU(N) using
bounded operators. {\it Invent. Math.} {\bf133} (1998) 467--538.

\bibitem{xu 2003} Xu F.: Strong additivity and conformal nets. 
math.QA/0303266. 

}

\end{thebibliography}
\end{document}